\documentclass{amsart}

\usepackage{graphicx}

\usepackage{amssymb}

\newtheorem{theorem}{Theorem} [section]

\newtheorem{thm}[theorem]{Theorem}

\newtheorem{cor}[theorem]{Corollary}

\newtheorem{definition}[theorem]{Definition}

\newtheorem{lemma}[theorem]{Lemma}

\newtheorem{prop}[theorem]{Proposition}

\newtheorem{obs}[theorem]{Observation}

\newtheorem{remark}[theorem]{Remark}

\begin{document}

%\date{\today}

\title{Classifying Intrinsically Linked Tournaments by Score Sequence}


\author{Thomas Fleming}
\author{Joel Foisy}

\begin{abstract}
A tournament on 8 or more vertices may be intrinsically linked as a directed graph.  We begin the classification of intrinsically linked tournaments by examining their score sequences.  While many distinct tournaments may have the same score sequence, there exist score sequences $S$ such that any tournament with score sequence $S$ has an embedding with no nonsplit consistently oriented link.   We call such score sequences \emph{linkless}, and we show that the vast majority of score sequences for 8 vertex tournaments are linkless.  

We also extend these results to $n$ vertex tournaments and are able to classify many longer score sequences as well. We show that for any $n$, there exist at least $O(n)$ linkless score sequences, but we conjecture that the fraction of score sequences of length $n$ that are linkless goes to 0 as $n$ becomes large.
\end{abstract}

\maketitle

\section{Introduction}

A graph is called \emph{intrinsically linked} if every embedding of that graph into $S^3$ contains disjoint cycles that form a non-split link.   This property was first studied in \cite{sachs} and \cite{cg}.  A directed graph $G$ is said to be \emph{intrinsically linked as a directed graph} if every embedding of $G$ into $S^3$ contains cycles that form a non-split link $L$, and the edges of $G$ that make up each component of $L$ have a consistent orientation.  Intrinsically linked directed graphs were first studied in \cite{FHR}.  The existence of intrinsically linked directed graphs with linking and knotting structures as complex as the unoriented graph case has been established \cite{mnp}.

In the case of undirected graphs, if $H$ is a minor of $G$ and $G$ has a linkless embedding, then $H$ does as well \cite{nt}.  Thus there is a finite family of minor minimal obstructions to a linkless embedding \cite{rs}, and this has been characterized as the Petersen family of graphs \cite{rst}.  Thus, given an undirected graph, we may determine if it is intrinsically linked by checking if it contains a Petersen family graph as a minor.   For example, any graph with $n \geq 6$ vertices is intrinsically linked if it contains more than $4n - 10$ edges, as this implies a $K_6$ minor \cite{mader}, and $K_6$ is in the Petersen family.  
 
For directed graphs, linkless embeddings are not preserved by the minor operation \cite{FHR}, but are preserved under certain other moves \cite{flemfoisy}.  Thus, not only is a classification of minor minimal examples for intrinsic linking in directed graphs unknown, it is unclear which minor-like operation (if any) will allow such a classification.  Given an arbitrary directed graph, it is currently difficult to determine whether it is intrinsically linked as a directed graph or not, except for the case of very dense graphs \cite{flemfoisy}.  However, we will show that for a \emph{tournament} it is often easy to determine that it is not intrinsically linked.

A tournament is a directed graph with exactly one directed edge between each pair of vertices.  Equivalently, a tournament is a choice of orientation for the edges of a complete graph $K_n$.  Tournaments may be intrinsically linked as a directed graphs \cite{flemfoisy2}, however, in contrast to undirected graphs, it is possible to have arbitrarily large tournaments that are not intrinsically linked as a directed graph.   For example, a \emph{transitive tournament} is a tournament where if edge $ab$ is oriented from $a$ to $b$ and edge $bc$ is oriented from $b$ to $c$, then edge $ac$ is oriented from $a$ to $c$.  In particular, a transitive tournament contains no consistently oriented cycles, and hence cannot be intrinsically linked as a directed graph no matter how many vertices it contains.  

The \emph{score sequence} of a tournament is the outdegree of the vertices, listed in non-decreasing order.  A transitive tournament on $n$ vertices can be identified by its score sequence, specifically $(0, 1, 2 \ldots n-2, n-1)$.   This motivates our approach:  given a score sequence for a tournament $T$, can we determine if $T$ has a linkless embedding?     

We say that a score sequence $S$ is \emph{linkless} if any tournament $T$ with score sequence $S$ has an embedding that contains no non-split consistently oriented link. We say that a score sequence $S'$ \emph{has an intrinsically linked representative} if there exists a tournament $T'$ with score sequence $S'$, where $T'$ is intrinsically linked as a directed graph.  

Any tournament on 7 or fewer vertices is not intrinsically linked as a directed graph, and there exists a tournament on 8 vertices that is intrinsically linked as a directed graph \cite{flemfoisy2}.  Thus 8 vertex tournaments are the first case where some score sequences are linkless, and some score sequences have intrinsically linked representatives.   We classify most score sequences for tournaments on 8 vertices in Section \ref{8_vert_section}, and extend the results to score sequences for tournaments with a larger number of vertices in Section \ref{9_vert_section}.  

The techniques of Sections \ref{8_vert_section} and \ref{9_vert_section} allow us to classify 162 of 167 score sequences for 8 vertex tournaments into those that are linkless (147 sequences) and those with an intrinsically linked representative (15 sequences). These results are summarized in Table \ref{8_vert_table}.  We may similarly classify 453 of 490 score sequences for 9 vertex tournaments, 1336 of 1486 for 10 vertex tournaments, and 4127 of 4639 for 11 vertex tournaments.  

We show in Section \ref{further_sec} that if a score sequence $S$ has an intrinsically linked representative, there can exist tournaments $T$ and $T'$ with score sequence $S$ where $T$ is intrinsically linked as a directed graph, and $T'$ has a linkless embedding.   Thus, if $S$ has an intrinsically linked representative, it does not imply that a tournament with score sequence $S$ is intrinsically linked as a directed graph.  We conjecture that for some $n$, there exists a score sequence $S$ such that if an $n$ vertex tournament $T$ has score sequence $S$, then $T$ is intrinsically linked as a directed graph.

In Section \ref{further_sec}, we also give a lower bound on the number of linkless score sequences for $n$ vertex tournaments.  We conjecture that the fraction of length $n$ score sequences that are linkless approaches 0 as $n$ goes to infinity.   That is, we conjecture that as $n$ becomes large, almost all score sequences have an intrinsically linked representative.

Landau's theorem \cite{landau} allows us to enumerate all score sequences for 8 vertex tournaments. Table \ref{8_vert_table} shows the classification of these sequences using our results.  Python code is available from the authors that will a) given a score sequence, report whether it is linkless, has an intrinsically linked representative, or that its intrinsic linking status is unknown,  b) given $n$, produce a list of all score sequences of length $n$ and their intrinsic linking status, and c) given $n$ and sequence fragments, return all score sequences of length $n$ that contain the fragments and their intrinsic linking status.

\begin{table}[]
\scalebox{0.6}{
\begin{tabular}{lllllllllll}
\cline{1-3} \cline{5-7} \cline{9-11}
\multicolumn{3}{|c|}{8 Vertex Tournaments} & \multicolumn{1}{l|}{} & \multicolumn{3}{c|}{8 Vertex Tournaments} & \multicolumn{1}{l|}{} & \multicolumn{3}{c|}{8 Vertex Tournaments} \\ \cline{1-3} \cline{5-7} \cline{9-11} 
\multicolumn{1}{|l|}{Score Sequence} & \multicolumn{1}{l|}{Status} & \multicolumn{1}{l|}{Reason} & \multicolumn{1}{l|}{} & \multicolumn{1}{l|}{Score Sequence} & \multicolumn{1}{l|}{Status} & \multicolumn{1}{l|}{Reason} & \multicolumn{1}{l|}{} & \multicolumn{1}{l|}{Score Sequence} & \multicolumn{1}{l|}{Status} & \multicolumn{1}{l|}{Reason} \\ \cline{1-3} \cline{5-7} \cline{9-11} 
\multicolumn{1}{|l|}{(0, 1, 2, 3, 4, 5, 6, 7) } & \multicolumn{1}{l|}{ linkless } & \multicolumn{1}{l|}{\ref{lemma_0_reduce}} & \multicolumn{1}{l|}{ } & \multicolumn{1}{l|}{(0, 3, 3, 4, 4, 4, 5, 5) } & \multicolumn{1}{l|}{ linkless } & \multicolumn{1}{l|}{\ref{lemma_0_reduce}} & \multicolumn{1}{l|}{ } & \multicolumn{1}{l|}{(1, 2, 3, 3, 3, 4, 6, 6) } & \multicolumn{1}{l|}{ linkless } & \multicolumn{1}{l|}{\ref{lemma_0_reduce}}\\ \cline{1-3} \cline{5-7} \cline{9-11}
\multicolumn{1}{|l|}{(0, 1, 2, 3, 4, 6, 6, 6) } & \multicolumn{1}{l|}{ linkless } & \multicolumn{1}{l|}{\ref{lemma_0_reduce}} & \multicolumn{1}{l|}{ } & \multicolumn{1}{l|}{(0, 3, 4, 4, 4, 4, 4, 5) } & \multicolumn{1}{l|}{ linkless } & \multicolumn{1}{l|}{\ref{lemma_0_reduce}} & \multicolumn{1}{l|}{ } & \multicolumn{1}{l|}{(1, 2, 3, 3, 3, 5, 5, 6) } & \multicolumn{1}{l|}{ linkless } & \multicolumn{1}{l|}{\ref{lemma_0_reduce}}\\ \cline{1-3} \cline{5-7} \cline{9-11}
\multicolumn{1}{|l|}{(0, 1, 2, 3, 5, 5, 5, 7) } & \multicolumn{1}{l|}{ linkless } & \multicolumn{1}{l|}{\ref{lemma_0_reduce}} & \multicolumn{1}{l|}{ } & \multicolumn{1}{l|}{(0, 4, 4, 4, 4, 4, 4, 4) } & \multicolumn{1}{l|}{ linkless } & \multicolumn{1}{l|}{\ref{lemma_0_reduce}} & \multicolumn{1}{l|}{ } & \multicolumn{1}{l|}{(1, 2, 3, 3, 4, 4, 4, 7) } & \multicolumn{1}{l|}{ linkless } & \multicolumn{1}{l|}{\ref{lemma_0_reduce}}\\ \cline{1-3} \cline{5-7} \cline{9-11}
\multicolumn{1}{|l|}{(0, 1, 2, 3, 5, 5, 6, 6) } & \multicolumn{1}{l|}{ linkless } & \multicolumn{1}{l|}{\ref{lemma_0_reduce}} & \multicolumn{1}{l|}{ } & \multicolumn{1}{l|}{(1, 1, 1, 3, 4, 5, 6, 7) } & \multicolumn{1}{l|}{ linkless } & \multicolumn{1}{l|}{\ref{lemma_0_reduce}} & \multicolumn{1}{l|}{ } & \multicolumn{1}{l|}{(1, 2, 3, 3, 4, 4, 5, 6) } & \multicolumn{1}{l|}{ linkless } & \multicolumn{1}{l|}{ \ref{prop_special_linkless}

 }\\ \cline{1-3} \cline{5-7} \cline{9-11}
\multicolumn{1}{|l|}{(0, 1, 2, 4, 4, 4, 6, 7) } & \multicolumn{1}{l|}{ linkless } & \multicolumn{1}{l|}{\ref{lemma_0_reduce}} & \multicolumn{1}{l|}{ } & \multicolumn{1}{l|}{(1, 1, 1, 3, 4, 6, 6, 6) } & \multicolumn{1}{l|}{ linkless } & \multicolumn{1}{l|}{\ref{lemma_0_reduce}} & \multicolumn{1}{l|}{ } & \multicolumn{1}{l|}{(1, 2, 3, 3, 4, 5, 5, 5) } & \multicolumn{1}{l|}{ linkless } & \multicolumn{1}{l|}{ \ref{prop_special_linkless}}\\ \cline{1-3} \cline{5-7} \cline{9-11}
\multicolumn{1}{|l|}{(0, 1, 2, 4, 4, 5, 5, 7) } & \multicolumn{1}{l|}{ linkless } & \multicolumn{1}{l|}{\ref{lemma_0_reduce}} & \multicolumn{1}{l|}{ } & \multicolumn{1}{l|}{(1, 1, 1, 3, 5, 5, 5, 7) } & \multicolumn{1}{l|}{ linkless } & \multicolumn{1}{l|}{\ref{lemma_0_reduce}} & \multicolumn{1}{l|}{ } & \multicolumn{1}{l|}{(1, 2, 3, 4, 4, 4, 4, 6) } & \multicolumn{1}{l|}{ IL rep } & \multicolumn{1}{l|}{\ref{flex_prop}}\\ \cline{1-3} \cline{5-7} \cline{9-11}
\multicolumn{1}{|l|}{(0, 1, 2, 4, 4, 5, 6, 6) } & \multicolumn{1}{l|}{ linkless } & \multicolumn{1}{l|}{\ref{lemma_0_reduce}} & \multicolumn{1}{l|}{ } & \multicolumn{1}{l|}{(1, 1, 1, 3, 5, 5, 6, 6) } & \multicolumn{1}{l|}{ linkless } & \multicolumn{1}{l|}{\ref{lemma_0_reduce}} & \multicolumn{1}{l|}{ } & \multicolumn{1}{l|}{(1, 2, 3, 4, 4, 4, 5, 5) } & \multicolumn{1}{l|}{ IL rep } & \multicolumn{1}{l|}{\ref{flex_prop}}\\ \cline{1-3} \cline{5-7} \cline{9-11}
\multicolumn{1}{|l|}{(0, 1, 2, 4, 5, 5, 5, 6) } & \multicolumn{1}{l|}{ linkless } & \multicolumn{1}{l|}{\ref{lemma_0_reduce}} & \multicolumn{1}{l|}{ } & \multicolumn{1}{l|}{(1, 1, 1, 4, 4, 4, 6, 7) } & \multicolumn{1}{l|}{ linkless } & \multicolumn{1}{l|}{\ref{lemma_0_reduce}} & \multicolumn{1}{l|}{ } & \multicolumn{1}{l|}{(1, 2, 4, 4, 4, 4, 4, 5) } & \multicolumn{1}{l|}{ linkless } & \multicolumn{1}{l|}{ \ref{prop_special_linkless}}\\ \cline{1-3} \cline{5-7} \cline{9-11}
\multicolumn{1}{|l|}{(0, 1, 2, 5, 5, 5, 5, 5) } & \multicolumn{1}{l|}{ linkless } & \multicolumn{1}{l|}{\ref{lemma_0_reduce}} & \multicolumn{1}{l|}{ } & \multicolumn{1}{l|}{(1, 1, 1, 4, 4, 5, 5, 7) } & \multicolumn{1}{l|}{ linkless } & \multicolumn{1}{l|}{\ref{lemma_0_reduce}} & \multicolumn{1}{l|}{ } & \multicolumn{1}{l|}{(1, 3, 3, 3, 3, 3, 5, 7) } & \multicolumn{1}{l|}{ linkless } & \multicolumn{1}{l|}{\ref{lemma_0_reduce}}\\ \cline{1-3} \cline{5-7} \cline{9-11}
\multicolumn{1}{|l|}{(0, 1, 3, 3, 3, 5, 6, 7) } & \multicolumn{1}{l|}{ linkless } & \multicolumn{1}{l|}{\ref{lemma_0_reduce}} & \multicolumn{1}{l|}{ } & \multicolumn{1}{l|}{(1, 1, 1, 4, 4, 5, 6, 6) } & \multicolumn{1}{l|}{ linkless } & \multicolumn{1}{l|}{\ref{lemma_0_reduce}} & \multicolumn{1}{l|}{ } & \multicolumn{1}{l|}{(1, 3, 3, 3, 3, 3, 6, 6) } & \multicolumn{1}{l|}{ linkless } & \multicolumn{1}{l|}{\ref{lemma_0_reduce}}\\ \cline{1-3} \cline{5-7} \cline{9-11}
\multicolumn{1}{|l|}{(0, 1, 3, 3, 3, 6, 6, 6) } & \multicolumn{1}{l|}{ linkless } & \multicolumn{1}{l|}{\ref{lemma_0_reduce}} & \multicolumn{1}{l|}{ } & \multicolumn{1}{l|}{(1, 1, 1, 4, 5, 5, 5, 6) } & \multicolumn{1}{l|}{ linkless } & \multicolumn{1}{l|}{\ref{lemma_0_reduce}} & \multicolumn{1}{l|}{ } & \multicolumn{1}{l|}{(1, 3, 3, 3, 3, 4, 4, 7) } & \multicolumn{1}{l|}{ linkless } & \multicolumn{1}{l|}{\ref{lemma_0_reduce}}\\ \cline{1-3} \cline{5-7} \cline{9-11}
\multicolumn{1}{|l|}{(0, 1, 3, 3, 4, 4, 6, 7) } & \multicolumn{1}{l|}{ linkless } & \multicolumn{1}{l|}{\ref{lemma_0_reduce}} & \multicolumn{1}{l|}{ } & \multicolumn{1}{l|}{(1, 1, 1, 5, 5, 5, 5, 5) } & \multicolumn{1}{l|}{ linkless } & \multicolumn{1}{l|}{\ref{lemma_0_reduce}} & \multicolumn{1}{l|}{ } & \multicolumn{1}{l|}{(1, 3, 3, 3, 3, 4, 5, 6) } & \multicolumn{1}{l|}{ IL rep } & \multicolumn{1}{l|}{\ref{flex_prop}}\\ \cline{1-3} \cline{5-7} \cline{9-11}
\multicolumn{1}{|l|}{(0, 1, 3, 3, 4, 5, 5, 7) } & \multicolumn{1}{l|}{ linkless } & \multicolumn{1}{l|}{\ref{lemma_0_reduce}} & \multicolumn{1}{l|}{ } & \multicolumn{1}{l|}{(1, 1, 2, 2, 4, 5, 6, 7) } & \multicolumn{1}{l|}{ linkless } & \multicolumn{1}{l|}{\ref{lemma_0_reduce}} & \multicolumn{1}{l|}{ } & \multicolumn{1}{l|}{(1, 3, 3, 3, 3, 5, 5, 5) } & \multicolumn{1}{l|}{ linkless } & \multicolumn{1}{l|}{ \ref{prop_special_linkless}}\\ \cline{1-3} \cline{5-7} \cline{9-11}
\multicolumn{1}{|l|}{(0, 1, 3, 3, 4, 5, 6, 6) } & \multicolumn{1}{l|}{ linkless } & \multicolumn{1}{l|}{\ref{lemma_0_reduce}} & \multicolumn{1}{l|}{ } & \multicolumn{1}{l|}{(1, 1, 2, 2, 4, 6, 6, 6) } & \multicolumn{1}{l|}{ linkless } & \multicolumn{1}{l|}{\ref{lemma_0_reduce}} & \multicolumn{1}{l|}{ } & \multicolumn{1}{l|}{(1, 3, 3, 3, 4, 4, 4, 6) } & \multicolumn{1}{l|}{ linkless } & \multicolumn{1}{l|}{\ref{prop_special_linkless} }\\ \cline{1-3} \cline{5-7} \cline{9-11}
\multicolumn{1}{|l|}{(0, 1, 3, 3, 5, 5, 5, 6) } & \multicolumn{1}{l|}{ linkless } & \multicolumn{1}{l|}{\ref{lemma_0_reduce}} & \multicolumn{1}{l|}{ } & \multicolumn{1}{l|}{(1, 1, 2, 2, 5, 5, 5, 7) } & \multicolumn{1}{l|}{ linkless } & \multicolumn{1}{l|}{\ref{lemma_0_reduce}} & \multicolumn{1}{l|}{ } & \multicolumn{1}{l|}{(1, 3, 3, 3, 4, 4, 5, 5) } & \multicolumn{1}{l|}{ IL rep } & \multicolumn{1}{l|}{\ref{flex_prop}}\\ \cline{1-3} \cline{5-7} \cline{9-11}
\multicolumn{1}{|l|}{(0, 1, 3, 4, 4, 4, 5, 7) } & \multicolumn{1}{l|}{ linkless } & \multicolumn{1}{l|}{\ref{lemma_0_reduce}} & \multicolumn{1}{l|}{ } & \multicolumn{1}{l|}{(1, 1, 2, 2, 5, 5, 6, 6) } & \multicolumn{1}{l|}{ linkless } & \multicolumn{1}{l|}{\ref{lemma_0_reduce}} & \multicolumn{1}{l|}{ } & \multicolumn{1}{l|}{(1, 3, 3, 4, 4, 4, 4, 5) } & \multicolumn{1}{l|}{ IL rep } & \multicolumn{1}{l|}{\ref{flex_prop}}\\ \cline{1-3} \cline{5-7} \cline{9-11}
\multicolumn{1}{|l|}{(0, 1, 3, 4, 4, 4, 6, 6) } & \multicolumn{1}{l|}{ linkless } & \multicolumn{1}{l|}{\ref{lemma_0_reduce}} & \multicolumn{1}{l|}{ } & \multicolumn{1}{l|}{(1, 1, 2, 3, 3, 5, 6, 7) } & \multicolumn{1}{l|}{ linkless } & \multicolumn{1}{l|}{\ref{lemma_0_reduce}} & \multicolumn{1}{l|}{ } & \multicolumn{1}{l|}{(1, 3, 4, 4, 4, 4, 4, 4) } & \multicolumn{1}{l|}{ linkless } & \multicolumn{1}{l|}{ \ref{prop_special_linkless} }\\ \cline{1-3} \cline{5-7} \cline{9-11}
\multicolumn{1}{|l|}{(0, 1, 3, 4, 4, 5, 5, 6) } & \multicolumn{1}{l|}{ linkless } & \multicolumn{1}{l|}{\ref{lemma_0_reduce}} & \multicolumn{1}{l|}{ } & \multicolumn{1}{l|}{(1, 1, 2, 3, 3, 6, 6, 6) } & \multicolumn{1}{l|}{ linkless } & \multicolumn{1}{l|}{\ref{lemma_0_reduce}} & \multicolumn{1}{l|}{ } & \multicolumn{1}{l|}{(2, 2, 2, 2, 2, 5, 6, 7) } & \multicolumn{1}{l|}{ linkless } & \multicolumn{1}{l|}{\ref{lemma_0_reduce}}\\ \cline{1-3} \cline{5-7} \cline{9-11}
\multicolumn{1}{|l|}{(0, 1, 3, 4, 5, 5, 5, 5) } & \multicolumn{1}{l|}{ linkless } & \multicolumn{1}{l|}{\ref{lemma_0_reduce}} & \multicolumn{1}{l|}{ } & \multicolumn{1}{l|}{(1, 1, 2, 3, 4, 4, 6, 7) } & \multicolumn{1}{l|}{ linkless } & \multicolumn{1}{l|}{\ref{lemma_0_reduce}} & \multicolumn{1}{l|}{ } & \multicolumn{1}{l|}{(2, 2, 2, 2, 2, 6, 6, 6) } & \multicolumn{1}{l|}{ linkless } & \multicolumn{1}{l|}{\ref{lemma_0_reduce}}\\ \cline{1-3} \cline{5-7} \cline{9-11}
\multicolumn{1}{|l|}{(0, 1, 4, 4, 4, 4, 4, 7) } & \multicolumn{1}{l|}{ linkless } & \multicolumn{1}{l|}{\ref{lemma_0_reduce}} & \multicolumn{1}{l|}{ } & \multicolumn{1}{l|}{(1, 1, 2, 3, 4, 5, 5, 7) } & \multicolumn{1}{l|}{ linkless } & \multicolumn{1}{l|}{\ref{lemma_0_reduce}} & \multicolumn{1}{l|}{ } & \multicolumn{1}{l|}{(2, 2, 2, 2, 3, 4, 6, 7) } & \multicolumn{1}{l|}{ linkless } & \multicolumn{1}{l|}{\ref{lemma_0_reduce}}\\ \cline{1-3} \cline{5-7} \cline{9-11}
\multicolumn{1}{|l|}{(0, 1, 4, 4, 4, 4, 5, 6) } & \multicolumn{1}{l|}{ linkless } & \multicolumn{1}{l|}{\ref{lemma_0_reduce}} & \multicolumn{1}{l|}{ } & \multicolumn{1}{l|}{(1, 1, 2, 3, 4, 5, 6, 6) } & \multicolumn{1}{l|}{ linkless } & \multicolumn{1}{l|}{\ref{lemma_0_reduce}} & \multicolumn{1}{l|}{ } & \multicolumn{1}{l|}{(2, 2, 2, 2, 3, 5, 5, 7) } & \multicolumn{1}{l|}{ linkless } & \multicolumn{1}{l|}{\ref{lemma_0_reduce}}\\ \cline{1-3} \cline{5-7} \cline{9-11}
\multicolumn{1}{|l|}{(0, 1, 4, 4, 4, 5, 5, 5) } & \multicolumn{1}{l|}{ linkless } & \multicolumn{1}{l|}{\ref{lemma_0_reduce}} & \multicolumn{1}{l|}{ } & \multicolumn{1}{l|}{(1, 1, 2, 3, 5, 5, 5, 6) } & \multicolumn{1}{l|}{ linkless } & \multicolumn{1}{l|}{\ref{lemma_0_reduce}} & \multicolumn{1}{l|}{ } & \multicolumn{1}{l|}{(2, 2, 2, 2, 3, 5, 6, 6) } & \multicolumn{1}{l|}{ linkless } & \multicolumn{1}{l|}{\ref{lemma_0_reduce}}\\ \cline{1-3} \cline{5-7} \cline{9-11}
\multicolumn{1}{|l|}{(0, 2, 2, 2, 4, 5, 6, 7) } & \multicolumn{1}{l|}{ linkless } & \multicolumn{1}{l|}{\ref{lemma_0_reduce}} & \multicolumn{1}{l|}{ } & \multicolumn{1}{l|}{(1, 1, 2, 4, 4, 4, 5, 7) } & \multicolumn{1}{l|}{ linkless } & \multicolumn{1}{l|}{\ref{lemma_0_reduce}} & \multicolumn{1}{l|}{ } & \multicolumn{1}{l|}{(2, 2, 2, 2, 4, 4, 5, 7) } & \multicolumn{1}{l|}{ linkless } & \multicolumn{1}{l|}{\ref{lemma_0_reduce}}\\ \cline{1-3} \cline{5-7} \cline{9-11}
\multicolumn{1}{|l|}{(0, 2, 2, 2, 4, 6, 6, 6) } & \multicolumn{1}{l|}{ linkless } & \multicolumn{1}{l|}{\ref{lemma_0_reduce}} & \multicolumn{1}{l|}{ } & \multicolumn{1}{l|}{(1, 1, 2, 4, 4, 4, 6, 6) } & \multicolumn{1}{l|}{ linkless } & \multicolumn{1}{l|}{\ref{lemma_0_reduce}} & \multicolumn{1}{l|}{ } & \multicolumn{1}{l|}{(2, 2, 2, 2, 4, 4, 6, 6) } & \multicolumn{1}{l|}{ linkless } & \multicolumn{1}{l|}{\ref{lemma_0_reduce}}\\ \cline{1-3} \cline{5-7} \cline{9-11}
\multicolumn{1}{|l|}{(0, 2, 2, 2, 5, 5, 5, 7) } & \multicolumn{1}{l|}{ linkless } & \multicolumn{1}{l|}{\ref{lemma_0_reduce}} & \multicolumn{1}{l|}{ } & \multicolumn{1}{l|}{(1, 1, 2, 4, 4, 5, 5, 6) } & \multicolumn{1}{l|}{ linkless } & \multicolumn{1}{l|}{\ref{lemma_0_reduce}} & \multicolumn{1}{l|}{ } & \multicolumn{1}{l|}{(2, 2, 2, 2, 4, 5, 5, 6) } & \multicolumn{1}{l|}{ linkless } & \multicolumn{1}{l|}{\ref{cor_8_out}}\\ \cline{1-3} \cline{5-7} \cline{9-11}
\multicolumn{1}{|l|}{(0, 2, 2, 2, 5, 5, 6, 6) } & \multicolumn{1}{l|}{ linkless } & \multicolumn{1}{l|}{\ref{lemma_0_reduce}} & \multicolumn{1}{l|}{ } & \multicolumn{1}{l|}{(1, 1, 2, 4, 5, 5, 5, 5) } & \multicolumn{1}{l|}{ linkless } & \multicolumn{1}{l|}{\ref{lemma_0_reduce}} & \multicolumn{1}{l|}{ } & \multicolumn{1}{l|}{(2, 2, 2, 2, 5, 5, 5, 5) } & \multicolumn{1}{l|}{ linkless } & \multicolumn{1}{l|}{\ref{cor_8_out}}\\ \cline{1-3} \cline{5-7} \cline{9-11}
\multicolumn{1}{|l|}{(0, 2, 2, 3, 3, 5, 6, 7) } & \multicolumn{1}{l|}{ linkless } & \multicolumn{1}{l|}{\ref{lemma_0_reduce}} & \multicolumn{1}{l|}{ } & \multicolumn{1}{l|}{(1, 1, 3, 3, 3, 4, 6, 7) } & \multicolumn{1}{l|}{ linkless } & \multicolumn{1}{l|}{\ref{lemma_0_reduce}} & \multicolumn{1}{l|}{ } & \multicolumn{1}{l|}{(2, 2, 2, 3, 3, 3, 6, 7) } & \multicolumn{1}{l|}{ linkless } & \multicolumn{1}{l|}{\ref{lemma_0_reduce}}\\ \cline{1-3} \cline{5-7} \cline{9-11}
\multicolumn{1}{|l|}{(0, 2, 2, 3, 3, 6, 6, 6) } & \multicolumn{1}{l|}{ linkless } & \multicolumn{1}{l|}{\ref{lemma_0_reduce}} & \multicolumn{1}{l|}{ } & \multicolumn{1}{l|}{(1, 1, 3, 3, 3, 5, 5, 7) } & \multicolumn{1}{l|}{ linkless } & \multicolumn{1}{l|}{\ref{lemma_0_reduce}} & \multicolumn{1}{l|}{ } & \multicolumn{1}{l|}{(2, 2, 2, 3, 3, 4, 5, 7) } & \multicolumn{1}{l|}{ linkless } & \multicolumn{1}{l|}{\ref{lemma_0_reduce}}\\ \cline{1-3} \cline{5-7} \cline{9-11}
\multicolumn{1}{|l|}{(0, 2, 2, 3, 4, 4, 6, 7) } & \multicolumn{1}{l|}{ linkless } & \multicolumn{1}{l|}{\ref{lemma_0_reduce}} & \multicolumn{1}{l|}{ } & \multicolumn{1}{l|}{(1, 1, 3, 3, 3, 5, 6, 6) } & \multicolumn{1}{l|}{ linkless } & \multicolumn{1}{l|}{\ref{lemma_0_reduce}} & \multicolumn{1}{l|}{ } & \multicolumn{1}{l|}{(2, 2, 2, 3, 3, 4, 6, 6) } & \multicolumn{1}{l|}{ linkless } & \multicolumn{1}{l|}{\ref{lemma_0_reduce}}\\ \cline{1-3} \cline{5-7} \cline{9-11}
\multicolumn{1}{|l|}{(0, 2, 2, 3, 4, 5, 5, 7) } & \multicolumn{1}{l|}{ linkless } & \multicolumn{1}{l|}{\ref{lemma_0_reduce}} & \multicolumn{1}{l|}{ } & \multicolumn{1}{l|}{(1, 1, 3, 3, 4, 4, 5, 7) } & \multicolumn{1}{l|}{ linkless } & \multicolumn{1}{l|}{\ref{lemma_0_reduce}} & \multicolumn{1}{l|}{ } & \multicolumn{1}{l|}{(2, 2, 2, 3, 3, 5, 5, 6) } & \multicolumn{1}{l|}{ linkless } & \multicolumn{1}{l|}{ \ref{prop_special_linkless} }\\ \cline{1-3} \cline{5-7} \cline{9-11}
\multicolumn{1}{|l|}{(0, 2, 2, 3, 4, 5, 6, 6) } & \multicolumn{1}{l|}{ linkless } & \multicolumn{1}{l|}{\ref{lemma_0_reduce}} & \multicolumn{1}{l|}{ } & \multicolumn{1}{l|}{(1, 1, 3, 3, 4, 4, 6, 6) } & \multicolumn{1}{l|}{ linkless } & \multicolumn{1}{l|}{\ref{lemma_0_reduce}} & \multicolumn{1}{l|}{ } & \multicolumn{1}{l|}{(2, 2, 2, 3, 4, 4, 4, 7) } & \multicolumn{1}{l|}{ linkless } & \multicolumn{1}{l|}{\ref{lemma_0_reduce}}\\ \cline{1-3} \cline{5-7} \cline{9-11}
\multicolumn{1}{|l|}{(0, 2, 2, 3, 5, 5, 5, 6) } & \multicolumn{1}{l|}{ linkless } & \multicolumn{1}{l|}{\ref{lemma_0_reduce}} & \multicolumn{1}{l|}{ } & \multicolumn{1}{l|}{(1, 1, 3, 3, 4, 5, 5, 6) } & \multicolumn{1}{l|}{ linkless } & \multicolumn{1}{l|}{\ref{lemma_0_reduce}} & \multicolumn{1}{l|}{ } & \multicolumn{1}{l|}{(2, 2, 2, 3, 4, 4, 5, 6) } & \multicolumn{1}{l|}{ linkless } & \multicolumn{1}{l|}{ \ref{prop_special_linkless} }\\ \cline{1-3} \cline{5-7} \cline{9-11}
\multicolumn{1}{|l|}{(0, 2, 2, 4, 4, 4, 5, 7) } & \multicolumn{1}{l|}{ linkless } & \multicolumn{1}{l|}{\ref{lemma_0_reduce}} & \multicolumn{1}{l|}{ } & \multicolumn{1}{l|}{(1, 1, 3, 3, 5, 5, 5, 5) } & \multicolumn{1}{l|}{ linkless } & \multicolumn{1}{l|}{\ref{lemma_0_reduce}} & \multicolumn{1}{l|}{ } & \multicolumn{1}{l|}{(2, 2, 2, 3, 4, 5, 5, 5) } & \multicolumn{1}{l|}{ unknown } & \multicolumn{1}{l|}{ }\\ \cline{1-3} \cline{5-7} \cline{9-11}
\multicolumn{1}{|l|}{(0, 2, 2, 4, 4, 4, 6, 6) } & \multicolumn{1}{l|}{ linkless } & \multicolumn{1}{l|}{\ref{lemma_0_reduce}} & \multicolumn{1}{l|}{ } & \multicolumn{1}{l|}{(1, 1, 3, 4, 4, 4, 4, 7) } & \multicolumn{1}{l|}{ linkless } & \multicolumn{1}{l|}{\ref{lemma_0_reduce}} & \multicolumn{1}{l|}{ } & \multicolumn{1}{l|}{(2, 2, 2, 4, 4, 4, 4, 6) } & \multicolumn{1}{l|}{ linkless } & \multicolumn{1}{l|}{ \ref{prop_special_linkless} }\\ \cline{1-3} \cline{5-7} \cline{9-11}
\multicolumn{1}{|l|}{(0, 2, 2, 4, 4, 5, 5, 6) } & \multicolumn{1}{l|}{ linkless } & \multicolumn{1}{l|}{\ref{lemma_0_reduce}} & \multicolumn{1}{l|}{ } & \multicolumn{1}{l|}{(1, 1, 3, 4, 4, 4, 5, 6) } & \multicolumn{1}{l|}{ linkless } & \multicolumn{1}{l|}{\ref{lemma_0_reduce}} & \multicolumn{1}{l|}{ } & \multicolumn{1}{l|}{(2, 2, 2, 4, 4, 4, 5, 5) } & \multicolumn{1}{l|}{ unknown } & \multicolumn{1}{l|}{ }\\ \cline{1-3} \cline{5-7} \cline{9-11}
\multicolumn{1}{|l|}{(0, 2, 2, 4, 5, 5, 5, 5) } & \multicolumn{1}{l|}{ linkless } & \multicolumn{1}{l|}{\ref{lemma_0_reduce}} & \multicolumn{1}{l|}{ } & \multicolumn{1}{l|}{(1, 1, 3, 4, 4, 5, 5, 5) } & \multicolumn{1}{l|}{ linkless } & \multicolumn{1}{l|}{\ref{lemma_0_reduce}} & \multicolumn{1}{l|}{ } & \multicolumn{1}{l|}{(2, 2, 3, 3, 3, 3, 5, 7) } & \multicolumn{1}{l|}{ linkless } & \multicolumn{1}{l|}{\ref{lemma_0_reduce}}\\ \cline{1-3} \cline{5-7} \cline{9-11}
\multicolumn{1}{|l|}{(0, 2, 3, 3, 3, 4, 6, 7) } & \multicolumn{1}{l|}{ linkless } & \multicolumn{1}{l|}{\ref{lemma_0_reduce}} & \multicolumn{1}{l|}{ } & \multicolumn{1}{l|}{(1, 1, 4, 4, 4, 4, 4, 6) } & \multicolumn{1}{l|}{ linkless } & \multicolumn{1}{l|}{\ref{lemma_0_reduce}} & \multicolumn{1}{l|}{ } & \multicolumn{1}{l|}{(2, 2, 3, 3, 3, 3, 6, 6) } & \multicolumn{1}{l|}{ linkless } & \multicolumn{1}{l|}{\ref{lemma_0_reduce}}\\ \cline{1-3} \cline{5-7} \cline{9-11}
\multicolumn{1}{|l|}{(0, 2, 3, 3, 3, 5, 5, 7) } & \multicolumn{1}{l|}{ linkless } & \multicolumn{1}{l|}{\ref{lemma_0_reduce}} & \multicolumn{1}{l|}{ } & \multicolumn{1}{l|}{(1, 1, 4, 4, 4, 4, 5, 5) } & \multicolumn{1}{l|}{ linkless } & \multicolumn{1}{l|}{\ref{lemma_0_reduce}} & \multicolumn{1}{l|}{ } & \multicolumn{1}{l|}{(2, 2, 3, 3, 3, 4, 4, 7) } & \multicolumn{1}{l|}{ linkless } & \multicolumn{1}{l|}{\ref{lemma_0_reduce}}\\ \cline{1-3} \cline{5-7} \cline{9-11}
\multicolumn{1}{|l|}{(0, 2, 3, 3, 3, 5, 6, 6) } & \multicolumn{1}{l|}{ linkless } & \multicolumn{1}{l|}{\ref{lemma_0_reduce}} & \multicolumn{1}{l|}{ } & \multicolumn{1}{l|}{(1, 2, 2, 2, 3, 5, 6, 7) } & \multicolumn{1}{l|}{ linkless } & \multicolumn{1}{l|}{\ref{lemma_0_reduce}} & \multicolumn{1}{l|}{ } & \multicolumn{1}{l|}{(2, 2, 3, 3, 3, 4, 5, 6) } & \multicolumn{1}{l|}{ IL rep } & \multicolumn{1}{l|}{\ref{flex_prop}}\\ \cline{1-3} \cline{5-7} \cline{9-11}
\multicolumn{1}{|l|}{(0, 2, 3, 3, 4, 4, 5, 7) } & \multicolumn{1}{l|}{ linkless } & \multicolumn{1}{l|}{\ref{lemma_0_reduce}} & \multicolumn{1}{l|}{ } & \multicolumn{1}{l|}{(1, 2, 2, 2, 3, 6, 6, 6) } & \multicolumn{1}{l|}{ linkless } & \multicolumn{1}{l|}{\ref{lemma_0_reduce}} & \multicolumn{1}{l|}{ } & \multicolumn{1}{l|}{(2, 2, 3, 3, 3, 5, 5, 5) } & \multicolumn{1}{l|}{ unknown } & \multicolumn{1}{l|}{ }\\ \cline{1-3} \cline{5-7} \cline{9-11}
\multicolumn{1}{|l|}{(0, 2, 3, 3, 4, 4, 6, 6) } & \multicolumn{1}{l|}{ linkless } & \multicolumn{1}{l|}{\ref{lemma_0_reduce}} & \multicolumn{1}{l|}{ } & \multicolumn{1}{l|}{(1, 2, 2, 2, 4, 4, 6, 7) } & \multicolumn{1}{l|}{ linkless } & \multicolumn{1}{l|}{\ref{lemma_0_reduce}} & \multicolumn{1}{l|}{ } & \multicolumn{1}{l|}{(2, 2, 3, 3, 4, 4, 4, 6) } & \multicolumn{1}{l|}{ IL rep } & \multicolumn{1}{l|}{\ref{flex_prop}}\\ \cline{1-3} \cline{5-7} \cline{9-11}
\multicolumn{1}{|l|}{(0, 2, 3, 3, 4, 5, 5, 6) } & \multicolumn{1}{l|}{ linkless } & \multicolumn{1}{l|}{\ref{lemma_0_reduce}} & \multicolumn{1}{l|}{ } & \multicolumn{1}{l|}{(1, 2, 2, 2, 4, 5, 5, 7) } & \multicolumn{1}{l|}{ linkless } & \multicolumn{1}{l|}{\ref{lemma_0_reduce}} & \multicolumn{1}{l|}{ } & \multicolumn{1}{l|}{(2, 2, 3, 3, 4, 4, 5, 5) } & \multicolumn{1}{l|}{ IL rep } & \multicolumn{1}{l|}{\ref{K332_prop}}\\ \cline{1-3} \cline{5-7} \cline{9-11}
\multicolumn{1}{|l|}{(0, 2, 3, 3, 5, 5, 5, 5) } & \multicolumn{1}{l|}{ linkless } & \multicolumn{1}{l|}{\ref{lemma_0_reduce}} & \multicolumn{1}{l|}{ } & \multicolumn{1}{l|}{(1, 2, 2, 2, 4, 5, 6, 6) } & \multicolumn{1}{l|}{ linkless } & \multicolumn{1}{l|}{\ref{lemma_0_reduce}} & \multicolumn{1}{l|}{ } & \multicolumn{1}{l|}{(2, 2, 3, 4, 4, 4, 4, 5) } & \multicolumn{1}{l|}{ IL rep } & \multicolumn{1}{l|}{\ref{K332_prop}}\\ \cline{1-3} \cline{5-7} \cline{9-11}
\multicolumn{1}{|l|}{(0, 2, 3, 4, 4, 4, 4, 7) } & \multicolumn{1}{l|}{ linkless } & \multicolumn{1}{l|}{\ref{lemma_0_reduce}} & \multicolumn{1}{l|}{ } & \multicolumn{1}{l|}{(1, 2, 2, 2, 5, 5, 5, 6) } & \multicolumn{1}{l|}{ linkless } & \multicolumn{1}{l|}{\ref{lemma_0_reduce}} & \multicolumn{1}{l|}{ } & \multicolumn{1}{l|}{(2, 2, 4, 4, 4, 4, 4, 4) } & \multicolumn{1}{l|}{ unknown } & \multicolumn{1}{l|}{ }\\ \cline{1-3} \cline{5-7} \cline{9-11}
\multicolumn{1}{|l|}{(0, 2, 3, 4, 4, 4, 5, 6) } & \multicolumn{1}{l|}{ linkless } & \multicolumn{1}{l|}{\ref{lemma_0_reduce}} & \multicolumn{1}{l|}{ } & \multicolumn{1}{l|}{(1, 2, 2, 3, 3, 4, 6, 7) } & \multicolumn{1}{l|}{ linkless } & \multicolumn{1}{l|}{\ref{lemma_0_reduce}} & \multicolumn{1}{l|}{ } & \multicolumn{1}{l|}{(2, 3, 3, 3, 3, 3, 4, 7) } & \multicolumn{1}{l|}{ linkless } & \multicolumn{1}{l|}{\ref{lemma_0_reduce}}\\ \cline{1-3} \cline{5-7} \cline{9-11}
\multicolumn{1}{|l|}{(0, 2, 3, 4, 4, 5, 5, 5) } & \multicolumn{1}{l|}{ linkless } & \multicolumn{1}{l|}{\ref{lemma_0_reduce}} & \multicolumn{1}{l|}{ } & \multicolumn{1}{l|}{(1, 2, 2, 3, 3, 5, 5, 7) } & \multicolumn{1}{l|}{ linkless } & \multicolumn{1}{l|}{\ref{lemma_0_reduce}} & \multicolumn{1}{l|}{ } & \multicolumn{1}{l|}{(2, 3, 3, 3, 3, 3, 5, 6) } & \multicolumn{1}{l|}{ linkless } & \multicolumn{1}{l|}{ \ref{prop_special_linkless}}\\ \cline{1-3} \cline{5-7} \cline{9-11}
\multicolumn{1}{|l|}{(0, 2, 4, 4, 4, 4, 4, 6) } & \multicolumn{1}{l|}{ linkless } & \multicolumn{1}{l|}{\ref{lemma_0_reduce}} & \multicolumn{1}{l|}{ } & \multicolumn{1}{l|}{(1, 2, 2, 3, 3, 5, 6, 6) } & \multicolumn{1}{l|}{ linkless } & \multicolumn{1}{l|}{\ref{lemma_0_reduce}} & \multicolumn{1}{l|}{ } & \multicolumn{1}{l|}{(2, 3, 3, 3, 3, 4, 4, 6) } & \multicolumn{1}{l|}{ IL rep } & \multicolumn{1}{l|}{\ref{flex_prop}}\\ \cline{1-3} \cline{5-7} \cline{9-11}
\multicolumn{1}{|l|}{(0, 2, 4, 4, 4, 4, 5, 5) } & \multicolumn{1}{l|}{ linkless } & \multicolumn{1}{l|}{\ref{lemma_0_reduce}} & \multicolumn{1}{l|}{ } & \multicolumn{1}{l|}{(1, 2, 2, 3, 4, 4, 5, 7) } & \multicolumn{1}{l|}{ linkless } & \multicolumn{1}{l|}{\ref{lemma_0_reduce}} & \multicolumn{1}{l|}{ } & \multicolumn{1}{l|}{(2, 3, 3, 3, 3, 4, 5, 5) } & \multicolumn{1}{l|}{ IL rep } & \multicolumn{1}{l|}{\ref{K332_prop}}\\ \cline{1-3} \cline{5-7} \cline{9-11}
\multicolumn{1}{|l|}{(0, 3, 3, 3, 3, 3, 6, 7) } & \multicolumn{1}{l|}{ linkless } & \multicolumn{1}{l|}{\ref{lemma_0_reduce}} & \multicolumn{1}{l|}{ } & \multicolumn{1}{l|}{(1, 2, 2, 3, 4, 4, 6, 6) } & \multicolumn{1}{l|}{ linkless } & \multicolumn{1}{l|}{\ref{lemma_0_reduce}} & \multicolumn{1}{l|}{ } & \multicolumn{1}{l|}{(2, 3, 3, 3, 4, 4, 4, 5) } & \multicolumn{1}{l|}{ IL rep } & \multicolumn{1}{l|}{\ref{K332_prop}}\\ \cline{1-3} \cline{5-7} \cline{9-11}
\multicolumn{1}{|l|}{(0, 3, 3, 3, 3, 4, 5, 7) } & \multicolumn{1}{l|}{ linkless } & \multicolumn{1}{l|}{\ref{lemma_0_reduce}} & \multicolumn{1}{l|}{ } & \multicolumn{1}{l|}{(1, 2, 2, 3, 4, 5, 5, 6) } & \multicolumn{1}{l|}{ linkless } & \multicolumn{1}{l|}{\ref{lemma_0_reduce}} & \multicolumn{1}{l|}{ } & \multicolumn{1}{l|}{(2, 3, 3, 4, 4, 4, 4, 4) } & \multicolumn{1}{l|}{ IL rep } & \multicolumn{1}{l|}{\ref{K332_prop}}\\ \cline{1-3} \cline{5-7} \cline{9-11}
\multicolumn{1}{|l|}{(0, 3, 3, 3, 3, 4, 6, 6) } & \multicolumn{1}{l|}{ linkless } & \multicolumn{1}{l|}{\ref{lemma_0_reduce}} & \multicolumn{1}{l|}{ } & \multicolumn{1}{l|}{(1, 2, 2, 3, 5, 5, 5, 5) } & \multicolumn{1}{l|}{ linkless } & \multicolumn{1}{l|}{\ref{cor_8_out}} & \multicolumn{1}{l|}{ } & \multicolumn{1}{l|}{(3, 3, 3, 3, 3, 3, 3, 7) } & \multicolumn{1}{l|}{ linkless } & \multicolumn{1}{l|}{\ref{lemma_0_reduce}}\\ \cline{1-3} \cline{5-7} \cline{9-11}
\multicolumn{1}{|l|}{(0, 3, 3, 3, 3, 5, 5, 6) } & \multicolumn{1}{l|}{ linkless } & \multicolumn{1}{l|}{\ref{lemma_0_reduce}} & \multicolumn{1}{l|}{ } & \multicolumn{1}{l|}{(1, 2, 2, 4, 4, 4, 4, 7) } & \multicolumn{1}{l|}{ linkless } & \multicolumn{1}{l|}{\ref{lemma_0_reduce}} & \multicolumn{1}{l|}{ } & \multicolumn{1}{l|}{(3, 3, 3, 3, 3, 3, 4, 6) } & \multicolumn{1}{l|}{ linkless } & \multicolumn{1}{l|}{ \ref{prop_special_linkless} }\\ \cline{1-3} \cline{5-7} \cline{9-11}
\multicolumn{1}{|l|}{(0, 3, 3, 3, 4, 4, 4, 7) } & \multicolumn{1}{l|}{ linkless } & \multicolumn{1}{l|}{\ref{lemma_0_reduce}} & \multicolumn{1}{l|}{ } & \multicolumn{1}{l|}{(1, 2, 2, 4, 4, 4, 5, 6) } & \multicolumn{1}{l|}{ linkless } & \multicolumn{1}{l|}{\ref{lemma_0_reduce}} & \multicolumn{1}{l|}{ } & \multicolumn{1}{l|}{(3, 3, 3, 3, 3, 3, 5, 5) } & \multicolumn{1}{l|}{ unknown } & \multicolumn{1}{l|}{ }\\ \cline{1-3} \cline{5-7} \cline{9-11}
\multicolumn{1}{|l|}{(0, 3, 3, 3, 4, 4, 5, 6) } & \multicolumn{1}{l|}{ linkless } & \multicolumn{1}{l|}{\ref{lemma_0_reduce}} & \multicolumn{1}{l|}{ } & \multicolumn{1}{l|}{(1, 2, 2, 4, 4, 5, 5, 5) } & \multicolumn{1}{l|}{ linkless } & \multicolumn{1}{l|}{\ref{prop_special_linkless} } & \multicolumn{1}{l|}{ } & \multicolumn{1}{l|}{(3, 3, 3, 3, 3, 4, 4, 5) } & \multicolumn{1}{l|}{ IL rep } & \multicolumn{1}{l|}{\ref{K332_prop}}\\ \cline{1-3} \cline{5-7} \cline{9-11}
\multicolumn{1}{|l|}{(0, 3, 3, 3, 4, 5, 5, 5) } & \multicolumn{1}{l|}{ linkless } & \multicolumn{1}{l|}{\ref{lemma_0_reduce}} & \multicolumn{1}{l|}{ } & \multicolumn{1}{l|}{(1, 2, 3, 3, 3, 3, 6, 7) } & \multicolumn{1}{l|}{ linkless } & \multicolumn{1}{l|}{\ref{lemma_0_reduce}} & \multicolumn{1}{l|}{ } & \multicolumn{1}{l|}{(3, 3, 3, 3, 4, 4, 4, 4) } & \multicolumn{1}{l|}{ IL rep } & \multicolumn{1}{l|}{\ref{K332_prop}}\\ \cline{1-3} \cline{5-7} \cline{9-11}
\multicolumn{1}{|l|}{(0, 3, 3, 4, 4, 4, 4, 6) } & \multicolumn{1}{l|}{ linkless } & \multicolumn{1}{l|}{\ref{lemma_0_reduce}} & \multicolumn{1}{l|}{ } & \multicolumn{1}{l|}{(1, 2, 3, 3, 3, 4, 5, 7) } & \multicolumn{1}{l|}{ linkless } & \multicolumn{1}{l|}{\ref{lemma_0_reduce}} & && & \\ \cline{1-3} \cline{5-7} 

 &  &  &  &  &  &  &  &  &  &  \\
 &  &  &  &  &  &  &  &  &  & 
\end{tabular}
}
\caption{The classification of score sequences for 8 vertex tournaments.} 
\label{8_vert_table}
\end{table}

\section{Preliminaries}

In this section we will introduce some of the techniques and results that will be useful in classifying tournaments using their score sequences.  First, notice that a tournament $T$ is intrinsically linked as a directed graph if and only if the tournament $\overline{T}$ is, where $\overline{T}$ is the tournament obtained by reversing the orientation of all edges of $T$.   Thus we need only classify a score sequence $S$ or the score sequence $\overline{S}$, which we will refer to as the \emph{dual score sequence} of $S$.

When studying a graph, it is often useful to look at minors of that graph.  Intrinsic linking in directed graphs is not well behaved under the standard graph minor operation \cite{FHR}.  In \cite{flemfoisy}, the authors introduce an operation called \emph{consistent edge contraction} that does preserve the property of having a linkless embedding.   We will make frequent use of consistent edge contraction in Sections \ref{8_vert_section} and \ref{9_vert_section}.

\begin{definition}
Let $e$ be an edge from $v$ to $w$ in a directed graph $G$ such that either $v$ is a sink in $G \setminus e$ or $w$ is a source in $G \setminus e$.  Let $H$ be the directed graph obtained from $G$ by deleting edge $e$ and forming a new vertex $v'$ by identifying $v$ and $w$.  We will say that $H$ is obtained from $G$ by \emph{consistent edge contraction}.
\end{definition}

\begin{figure}
\includegraphics[scale=0.65]{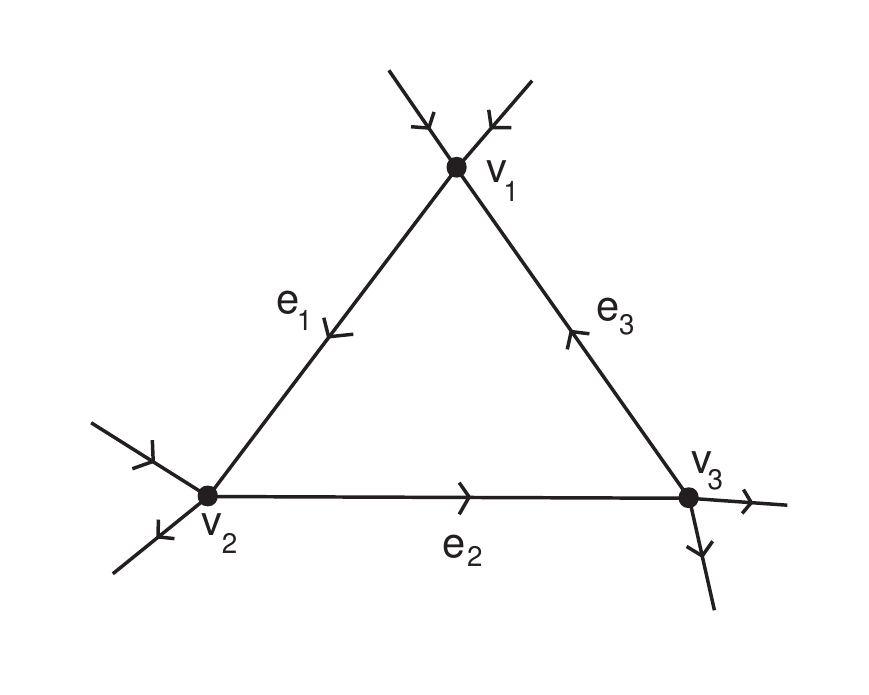}
\caption{Vertices $v_1$ and $v_2$ may be identified by consistent edge contraction, as $v_1$ is a sink in $G \setminus e_1$.  Vertices $v_2$ and $v_3$ may be identified by consistent edge contraction, as $v_3$ is a source in $G \setminus e_2$.  Vertices $v_1$ and $v_3$ cannot be identified by consistent edge contraction.} 
\label{edge_contraction_figure}
\end{figure}

It was shown in \cite{flemfoisy} that if $G$ is a directed graph that has an embedding with no nonsplit consistently oriented link, and $H$ is obtained from $G$ by a consistent edge contraction, then $H$ has an embedding with no nonsplit consistently oriented link as well.  In fact the relationship is much stronger, as shown below. 

\begin{thm}
If $H$ is obtained from $G$ by consistent edge contraction, then $H$
is intrinsically linked as a directed graph if and only if $G$ is. Further, $H$ is intrinsically knotted as a directed graph if and only if $G$ is.
\label{consist_edge_thm}
\end{thm}

\begin{proof}
Suppose $H$ is obtained from $G$ by consistent edge contraction on edge $e$ that
runs from vertex $v$ to vertex $w$.   We may assume that $w$ is a source in $G \setminus e$,
as the arguments below will proceed similarly if $v$ is a sink in $G \setminus e$.  We will
abuse notation and call the resulting vertex $w$ in $H$.

Upon contracting edge $e$, multiple edges with the same orientation or a loop may occur.  There are no edges of the form $xw$ in $G$ except for $vw$. Thus a multiple edge will occur when $wx$ and $vx$
are both edges of $G$. In this case, we delete the edge $vx$ from each pair.  A loop may be formed if there is an edge $e'$ from $w$
to $v$.  We delete this edge as well. Call the graph formed by deleting this set $E$ of edges $G'$. When
$e$ is contracted in $G'$, we obtain $H$ (with no loops or multiple edges).

We first consider intrinsic knotting. Let $\Gamma(G)$ denote the set of all consistently oriented cycles in $G$.

Suppose that there exists an embedding $f$ of $G$ that contains no consistently oriented
cycle that forms a nontrivial knot.  Note that $f$ is a knotless
embedding of $G'$ as $\Gamma(G') \subset \Gamma(G)$, so $f(G')$
contains no consistently oriented non-trivial knot as well. We may construct an embedding of $H$
by contracting edge $e$ until it and $v$ lie in a neighborhood of $w$ disjoint from the rest of $G'$.
Using this embedding, any cycle $c \in \Gamma(H)$ can be seen to be isotopic to
$c' \in \Gamma(G') \subset \Gamma(G)$.  Hence all such $c$ are trivial knots, and
$H$ has a knotless embedding if $G$ does.

Suppose that there exists an embedding $f$ of $H$ that contains no consistently oriented cycle that is
as non-trivial knot.   Then we may extend $f$ to an embedding of $G'$ by embedding edge $e$
and vertex $v$ in a neighborhood of vertex $w$. Any cycle $c \in \Gamma(G')$ is isotopic to a
cycle $c'$ in $\Gamma(H)$, and hence is a trivial knot.   We will now extend $f$ to an embedding
of $G$ so that any element of $\Gamma(G) \setminus \Gamma(G')$ is
isotopic to an element of $\Gamma(G')$ and hence trivial.  Note that a cycle $c \in
\Gamma(G) \setminus \Gamma(G')$ must contain an edge from $E$, and in fact can
contain at most one, as the vertex $w$ is a source, and $E$ contains only edges of the form $e'$ and $vx$.   A cycle with two edges of the form
$vx$ is not consistently oriented, and similarly a cycle that contains $e'$ and $vx$ is inconsistently oriented as well.

If $e'$ is in $E$, we may embed it in a neighborhood of
$e$, so that the cycle $e'e$ bounds a disk.  As $w$ is a source in $G \setminus e$, if a cycle $c \in \Gamma(G)$
contains $e'$ then $c=e'e$.  Suppose $vx \in E$ for some vertex $x$.  For a consistently oriented
cycle $c$, if $vx \in c$ then $w \notin c$.  Thus, there is some $c' \in \Gamma(G')$ where
$c = vxp$ and $c' = vwxp$ for some path $p$ in $G'$.  We may embed the edge $vx$ in a neighborhood
of the path $vwx$ so that $c$ is isotopic to $c'$ for any $p$.  Thus, any element of $\Gamma(G)$ is isotopic to
an element of $\Gamma(H)$ and so $G$ has a knotless embedding if $H$ does.

We now consider intrinsic linking.  Let $\Gamma^2(G)$ denote the set of all pairs of disjoint consistently oriented cycles in $G$.

Suppose that $G$ has an embedding $f$ where all elements of $\Gamma^2(G)$
are trivial links.  We may then delete the set of edges $E$ to form $f(G')$.  Clearly all elements
of $\Gamma^2(G')$ are trivial links as well.   We may then form an embedding $f$ of $H$
by contracting edge $e$ until it lies within a neighborhood of vertex $w$, disjoint from the rest of $G'$.
This gives an embedding of $H$, where any element of $\Gamma^2(H)$ is isotopic to an element in
$\Gamma^2(G')$, and hence a trivial link.   Thus $H$ has a linkless embedding if $G$ does.

Suppose that $H$ has an embedding $f$ where all elements of $\Gamma^2(H)$ are
trivial links.  We may form an embedding $f$ of $G'$ from $f(H)$ by embedding
$e$ within a neighborhood of $w$, disjoint from the rest of $H$.   Any element of
$\Gamma^2(G')$ is isotopic to an element of $\Gamma^2(H)$, and hence a trivial link.
We will now show that $f(G')$ can be extended to an embedding $f(G)$ so that
any element of $\Gamma^2(G)$ is isotopic to an element of $\Gamma^2(G')$, and hence trivial.

As $w$ is a source in $G \setminus e$, if an element $(L_1, L_2)$ of $\Gamma^2(G)$
contains both $v$ and $w$, $v \in L_1$ and $w \in L_1$.  If $e'$ is in $E$, then we may extend $f$ to $e$
by embedding $e'$ in a neighborhood of $e$ so that $e'e$ bounds a disk.   As before, since $w$ is a source
in $G \setminus e$, the only element of $\Gamma(G)$ that contains $e'$ is $e'e$.
Thus, if $e' \in (L_1,L_2) \in \Gamma^2(G)$, then $L_1 = e'e$ and the link is trivial.  Suppose
$vx \in E$.  We may extend $f$ by embedding the edge $vx$ in a neighborhood of the path $vwx$.  Then as before,
$c=vxp \in \Gamma(G)$ is isotopic to $c' = vwxp \in \Gamma(G')$.
If $(c,L_2) \in \Gamma^2(G)$, then $w \notin L_1, L_2$, and $(c, L_2) \cong (c',L_2) \in \Gamma^2(G')$.
Hence $(c, L_2)$ is a trivial link. Thus $f(G')$ can be extended to an embedding of $f(G)$ such that
any element of $\Gamma^2(G)$ is isotopic to an element of $\Gamma^2(H)$, and hence trivial.
Thus, $G$ has a linkless embedding if $H$ does.

\end{proof}

The following result is often useful for replacing an inconsistent cycle in a link with a consistent one.

\begin{figure}
\includegraphics[scale=0.15]{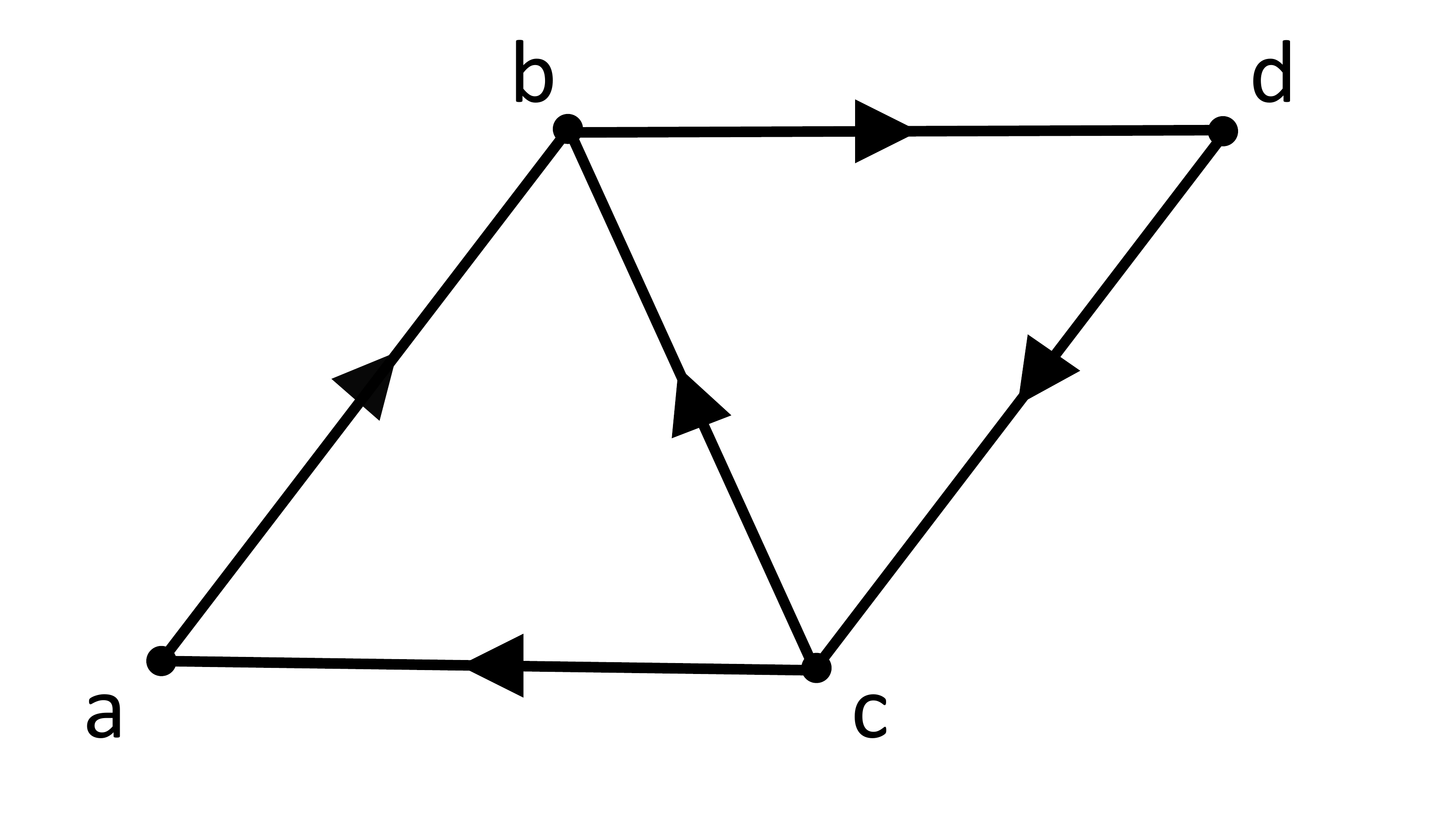}
\caption{Converting a linked inconsistent 3-cycle to a linked consistent cycle.} 
\label{fixing_fig}
\end{figure}

\begin{lemma}
Let $P_1, P_2$ and $P_3$ be consistently oriented paths in $G$, disjoint except for their endpoints $a$ and $b$.  Let $P_1$ and $P_2$ be oriented from $a$ to $b$, and $P_3$ be oriented from $b$ to $a$.  Let $C_1 = P_1 \cup P_2,  C_2 = P_2 \cup P_3$ and $C_3 = P_1 \cup P_3$.   If there exists a cycle $X$ disjoint from $P_i$ with $lk(X, C_1) \neq 0$, then $lk(X, C_2) \neq 0$ or $lk(X, C_3) \neq 0$.   Note that $C_2$ and $C_3$ are consistently oriented.

In particular, if $C_1 = abc$, $C_2 = bdc$, and $C_3 = abdc$ as shown in Figure \ref{fixing_fig}, and $X$ a cycle disjoint from each $C_i$ with $lk(X, C_1) \neq 0$, then $lk(X, C_2) \neq 0$ or $lk(X, C_3) \neq 0$ and $C_2, C_3$ are consistently oriented. 
\label{fix_orient_cor}
\end{lemma}

\begin{proof}
Considering $C_1$ as $P_1 \cup -P_2$ as an element of $H_1({\mathbb R}^3-X, {\mathbb Z})$, we have $ [C_1] + [C_2] - [C_3] = 0 \in H_1({\mathbb R}^3-X, {\mathbb Z})$.   As $[C_1]$ is not zero, one of $[C_2]$ and $[C_3]$ must be non-zero as well.
\end{proof}

Using consistent edge contraction, a tournament on 8 vertices can be reduced to a directed graph on 7 vertices, that will usually contain \emph{symmetric edges}.  We will call two edges between $v$ and $w$ symmetric edges if one is oriented from $v$ to $w$ and the other is oriented from $w$ to $v$.   If a directed graph $G$ has two edges between $v$ and $w$, with both directed from $v$ to $w$ or from $w$ to $v$, then $G$ is intrinsically linked as a directed graph if and only if $G'$ is, where $G'$ has a single edge between $v$ and $w$ and is otherwise identical to $G$, as a consistent cycle may use at most one of the edges between $v$ and $w$ in $G$.  If the consistent edge contraction of an 8 vertex tournament results in a directed graph $G$ that has multiple edges (but no symmetric edges) the graph $G$ can be reduced to a 7 vertex tournament and has a linkless embedding, as shown in \cite{flemfoisy2}.

\begin{thm}[Theorem 2.2 of \cite{flemfoisy2}] 
No tournament on 7 vertices is intrinsically linked as a directed graph.
\label{7_tourneys}
\end{thm}

When the 7 vertex directed graph has symmetric edges, we may choose embeddings where these pairs of symmetric edges bound disks, so it is often useful to understand the linked cycles in an embedding of $K_7$.

\begin{figure}
\includegraphics[scale=0.4]{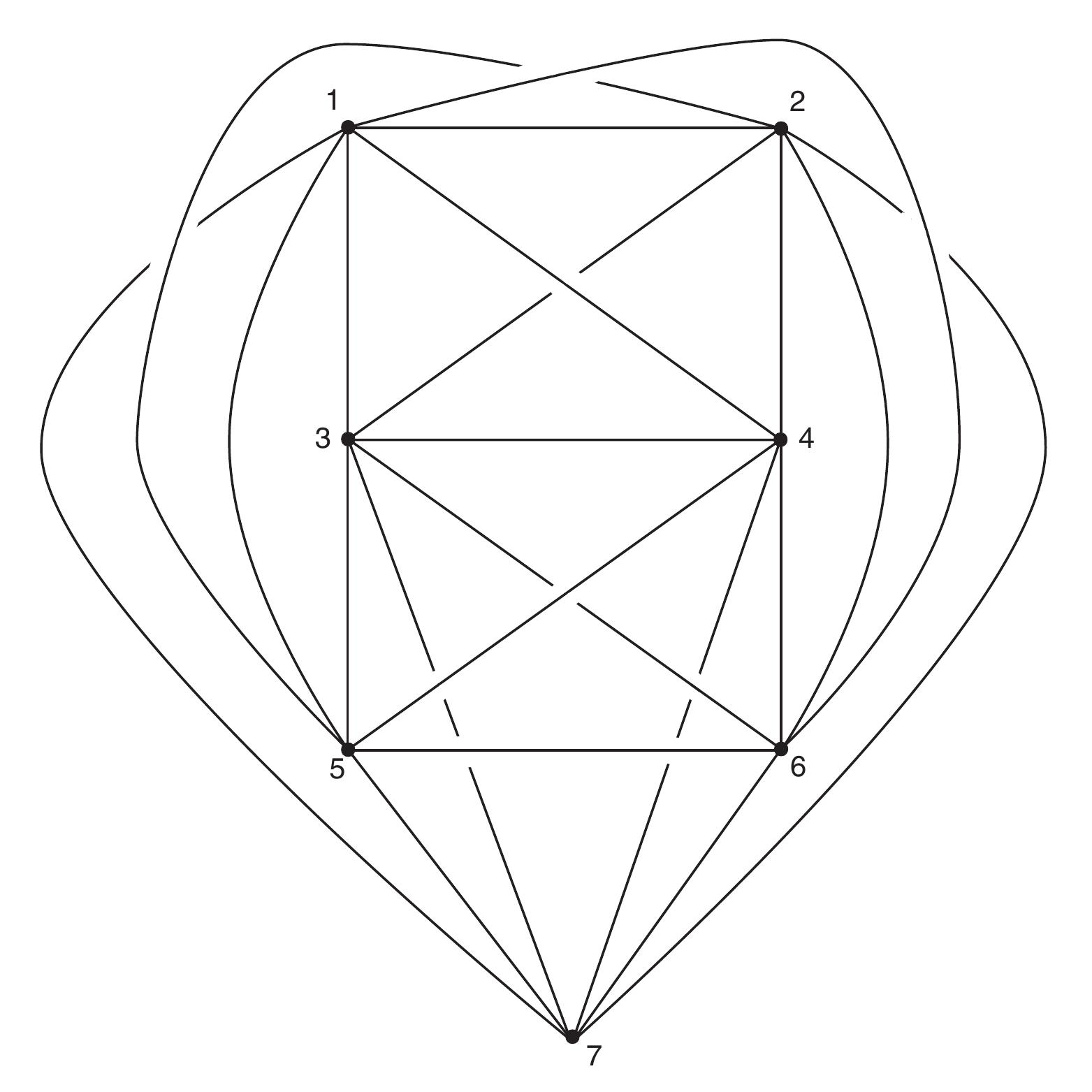}
\caption{The embedding of $K_7$ from \cite{cg}, which has exactly 21 non-split links \cite{fmellor}.} 
\label{FMellorK7}
\end{figure}

Figure \ref{FMellorK7} shows the embedding of $K_7$ from \cite{cg}, containing exactly 21 non-split links, which is the minimal number for any embedding of $K_7$ \cite{fmellor}.  We will refer to this as the CG embedding of $K_7$.   The linked cycles in the CG embedding are:

457-236 457-136 457-1362 457-1236

147-236 147-235 147-2356 147-2365

167-235 167-245 167-2435 167-2345

136-245 136-2547 136-2457

235-1467 235-1647

245-1376 245-1736

236-1475 236-1547

\begin{obs}
\label{LS_obs}
Every link in the CG embedding contains at least one of the cycles in the set $LS = \{ 236, 235, 136, 245, 1362, 1236, 2356, 2365, 2435, 2345\}$.    Thus, if each of these cycles has an inconsistent orientation, the CG embedding contains no consistently oriented non-split links.
\end{obs}

\medskip

A directed graph $G$ obtained from a consistent edge contraction on an 8 vertex tournament may have symmetric edges, but all such edges will be incident on a single vertex.   We may use the CG embedding to show that many of these graphs are not intrinsically linked as directed graphs by examining the score sequence of the 6 vertex tournament formed by deleting the preferred vertex.  We will use these lemmas in Section \ref{8_vert_section}.

\begin{lemma}
Let $G$ be a directed graph with 7 vertices, possibly with symmetric directed edges.  If there exists a vertex $v$ such that $G \setminus v$ is a tournament $T$ on 6 vertices, and the score sequence $S$ of $T$ contains $(0 \ldots)$ or $(\ldots 5)$ or $(1,1 \ldots)$ or $(\ldots 4,4)$, then $G$ is not intrinsically linked as a directed graph.
\label{sub_tourney_0_11}
\end{lemma}

\begin{proof}
We will label vertex $v$ as 7, choose vertex labels for $T$ and then place $G$ in the CG embedding (with any symmetric 
edges bounding disks). We will choose the vertex labels so that any non-split link in the CG embedding has at least one component 
that is not consistently oriented. By Observation \ref{LS_obs}, we need only examine cycles in the set $LS$.

Suppose $S$ contains 0 or 5. We need only consider the case of $0 \in S$, as the same proof applies to $5 \in S$ by reversing the orientation 
of every edge.  Label the vertex of out degree 0 as 2.  Then all elements of $LS$ have an inconsistent orientation except $136$.  As $T \setminus 2$
is a tournament on 5 vertices, there must be an inconsistently oriented 3-cycle.  Label its vertices 1, 3 and 6.  Thus $G$ may be placed in the CG
embedding with no non-split consistently oriented links.

Suppose $S$ contains 1,1 or 4,4.  As before, by reversing orientations if necessary, we need only consider 1,1. There is an edge between the vertices
of out degree 1.  Label the vertices so that this edge is oriented from 2 to 4. Then any cycle that includes 2, but not adjacently to 4, is inconsistently
oriented.  Thus, the only possible consistent elements of $LS$ are 136, 245, and 2435. Vertex 4 has an edge oriented to exactly 1 other vertex. 
Label that vertex 6. This forces 245 and 2435 to be inconsistent.  The vertices 1,3,5,6 form a tournament on four vertices, and vertex 6 is contained
in three 3-cycles, $6ab$, $6bc$ and $6ca$.  At least one of these must be inconsistent, say $6ab$. Label $a$ as 1 and $b$ as 3. Then all cycles in 
$LS$ are inconsistently oriented, and hence $G$ is not intrinsically linked as a directed graph.  
\end{proof}

\begin{lemma}
Let $G$ be a directed graph with 7 vertices, possibly with symmetric directed edges.  If there exists a vertex $v$ such that $G \setminus v$ is
a tournament $T$ on 6 vertices labeled $\{a,b,c,d,e,f\}$ whose score sequence is $(1,2,3,3,3,3)$, then either $G$ is not intrinsically linked as a 
directed graph, or $T$ has edge orientations $ba, ac, da, ea, fa, bc, db, eb, fb, cd, ce, cf, de, fd, ef$.

If $T$ has these edge orientations, and there is a single edge between $v$ and $b$, or a single edge between $v$ and $d$ with orientation $dv$, or a single edge between $v$ and each of $c$ and $d$, such that the edge orientations are $vc, vd$, then $G$ is not intrinsically linked as a 
directed graph.
\label{sub_tourney_123333}
\end{lemma}

\begin{proof}

We will label vertex $v$ as 7, choose vertex labels for $T$ and then place $G$ in the CG embedding (with any symmetric 
edges bounding disks). We will choose the vertex labels so that any non-split link in the CG embedding has at least one component 
that is not consistently oriented. By Observation \ref{LS_obs}, we need only examine cycles in the set $LS$.

The tournament $T$ has score sequence $(1,2,3,3,3,3)$. Label the vertex of out degree 1 as $a$, 
the vertex of out degree 2 as $b$ and the vertices of out degree 3 as $c,d,e,f$.  The vertices
$a$ and $b$ form a $K_2$, so there are 2 edges directed from $\{a,b\}$ to $\{c,d,e,f\}$.  There are two possible arrangements of $a,b$ in $T$.
First, edge $ab$ oriented from $a$ to $b$, and two edges from $b$ oriented to $\{c,d,e,f\}$.  Second, edge $ab$ oriented from $b$ to $a$, and 
one edge from each of $a$ and $b$ oriented to $\{c,d,e,f\}$. 

Vertices $c,d,e,f$ form a tournament $T'$ on 4 vertices, and each must have out degree 3 when including the edges to
$a$ and $b$.  Thus, no vertex may have out degree less than 1 or more than $3$ in $T'$.  Thus $T'$ must have score sequence $(1,1,1,3)$ or 
$(1,1,2,2)$.  Up to symmetry, there is a unique tournament with each score sequence.  Note that we cannot have edge $ab$ oriented from $a$ to $b$ 
and $T'$ with score sequence $(1,1,1,3)$, as both edges from $b$ to $\{c,d,e,f\}$ would have the same end points.  Thus, there are three cases 
to check. Case 1: edge $ab$ oriented from $a$ to $b$ and $T'$ with score sequence $(1,1,2,2)$.  Case 2: edge $ab$ oriented from $b$ to $a$ and $T'$ with score sequence $(1,1,2,2)$. Case 3: edge $ab$ oriented from $b$ to $a$ and $T'$ with score sequence $(1,1,1,3)$. 

In the first two cases, it is straightforward to choose vertex labels so that all of the cycles of $LS$ have inconsistent orientations.  Case 3 gives rise to the specific edge orientations of the hypothesis.

In this case no matter the choice of vertex labels, one or more elements of $LS$ have a consistent orientation.  Using the additional edge orientations between $v$ and vertices of $T$ required by the hypothesis, we may find a choice of vertex labels so that any consistently oriented cycle in $LS$ forms nonsplit links only with inconsistently oriented cycles. 

Thus, $G$ is not intrinsically linked as a directed graph.

\end{proof}

\begin{lemma}
Let $G$ be a directed graph with 7 vertices, possibly with symmetric directed edges.  If there exists a vertex $v$ such that $G \setminus v$ is
a tournament $T$ on 6 vertices whose score sequence is $(2,2,2,3,3,3)$ or $(1,2,2,3,3,4)$,  then $G$ is not intrinsically linked as a directed graph.
\label{sub_tourney_222333}
%\label{sub_tourney_122334}
\end{lemma}

The proof of Lemma \ref{sub_tourney_222333} is the similar to that of Lemma \ref{sub_tourney_0_11} and Lemma \ref{sub_tourney_123333}.

\medskip
In some cases we are able to reduce to a directed graph on 6 vertices.   In these cases the following consequences of \cite{FHR} can be useful.

The \emph{complete symmetric digraph} on $n$ vertices is the directed graph such that any pair of vertices $v$ and $w$ are joined by edges $vw$ and $wv$.  

\begin{thm}[Theorem 3.9 of \cite{FHR}]
Let $G$ be a directed graph formed by deleting $e_1$ and $e_2$ from the complete symmetric digraph on 6 vertices.  If $e_1$ and $e_2$ have a common vertex $v$ and are not consistently oriented, then $G$ is not intrinsically linked as a directed graph.
\label{FHR_edge_pair}
\end{thm}

\begin{cor}[Corollary 3.10 of \cite{FHR}]
A directed graph on 6 vertices with 23 or fewer edges is not intrinsically linked as a directed graph.
\label{FHR_23_edges}
\end{cor}

\section{Tournaments on 8 Vertices}
\label{8_vert_section}

In this section, we will classify score sequences for 8 vertex tournaments.  We know that some tournaments on 8 vertices are intrinsically linked \cite{flemfoisy2}.   We first pursue various constructions of intrinsically linked 8 vertex tournaments to find score sequences that have intrinsically linked representatives.

\begin{prop} The following score sequences and their dual score sequences have intrinsically linked representatives: $(2,3,3,4,4,4,4,4), (2,3,3,3,4,4,4,5)$, \linebreak $(3,3,3,3,4,4,4,4),$ $(3,3,3,3,3,4,4,5)$, $(2,2,3,3,4,4,5,5)$ and
$(2,2,3,4,4,4,4,5)$.
\label{K332_prop}
\end{prop}

\begin{proof}
We will construct intrinsically linked tournaments with the desired score sequences by iteratively assigning orientations to the edges of $K_8$.

Label the vertices of $K_8$ as $\{a_1, a_2, a_3, b_1, b_2, b_3, c_1, c_2\}$.  Consider the subgraph $H$ of $K_8$
isomorphic to $K_{3,3,2}$ formed by choosing the vertex partitions $A=\{a_1, a_2, a_3\}, B=\{b_1, b_2, b_3\}$ and
$C=\{c_1, c_2\}$.  Orient the edges of $H$ as follows:  from $A$ to $B$, from $B$ to $C$, from $C$ to $A$.

Every embedding of $K_{3,3,2}$ contains a pair of disjoint 3-cycles that have non-zero linking number
\cite{REU}.  As these 3-cycles are disjoint, they must be of the form $c_1a_ib_j$ and $c_2a_kb_l$.  By the
construction of $H$, these cycles are consistently oriented.  Thus, $H$ is intrinsically linked as a directed graph, and any tournament that contains $H$ as a subgraph will be intrinsically linked as a directed graph as well.

The vertices of $H$ have out degrees $(2,2,2,3,3,3,3,3)$.  We will now add edges to $H$ to construct intrinsically linked tournaments.  We may assume that edge $c_1c_2$ is oriented from $c_1$ to $c_2$, as otherwise we may exchange the labels of $c_1$ and $c_2$.  This gives out degrees $(2,2,2,3,3,3,3,4)$

We now need to assign orientations to the 3-cycles $a_1a_2a_3$ and $b_1b_2b_3$.  As all of the $a_i$ have out degree 3 and all of the $b_i$ have out degree 2, up to symmetry there are only 2 choices for each 3-cycle:  a consistent orientation or an inconsistent orientation.  In the first case, we add 1 out degree to each vertex, in the other we add 0, 1, and 2 respectively.  

We have four choices of orientations to complete the construction of $T$.  The cycle $a_1a_2a_3$ consistent and $b_1b_2b_3$ consistent, or $a_1a_2a_3$ consistent and $b_1b_2b_3$ inconsistent and so on.   These four options give rise to the first four score sequences.

To obtain the final two score sequences, start with $H$ and form a new directed graph $H'$ by reversing the orientations of the edges between $b_2$ and $C$ so that they are oriented from $c_i$ to $b_2$.

We can see that $H'$ is intrinsically linked as follows.  Find a pair of 3-cycles $L_1, L_2$ with non-zero linking number as before.   If these 3-cycles are both consistently oriented, we are done.

If a 3-cycle in the link is not consistently oriented, it must contain $b_2$.  Say it is $L_2$.
Orient the edge from $b_2$ to $b_1$ and from $b_2$ to $b_3$.  We may then assume that $L_1$ is $c_1a_1 b_1$
and $L_2$ is $c_2a_2 b_2$.  Form $L_3 = c_2a_2b_2b_3$ and $L_4 = c_2b_2b_3$.
By Lemma \ref{fix_orient_cor}, $L_1$ must have non-zero linking number with one of $L_3$ and $L_4$.  As these are both consistently
oriented, $H'$ is intrinsically linked as a directed graph.

Extending $H'$ to a tournament by choosing the orientation of the edges between the $a_i$ to form a consistent or inconsistent 3-cycle gives the two remaining score sequences.
\end{proof}

\begin{prop}
The score sequences $(2,3,3,3,3,4,4,6)$, $(1,3,3,3,3,4,5,6)$, \linebreak $(1,3,3,3,4,4,5,5)$, $(2,2,3,3,3,4,5,6)$ and their duals have intrinsically linked representatives.
\label{flex_prop}
% \label{ba_prop}

\end{prop}

\begin{figure}
\vskip -1in
\includegraphics [scale=.35]{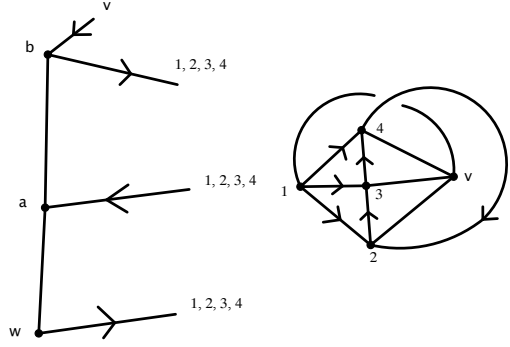}
\vskip -1in
\caption{A schematic of the directed graph $G$.}
\label{examp}
\end{figure}

\begin{proof}

We first consider the score sequence $(2,3,3,3,3,4,4,6)$.
Let $G$ be the directed graph depicted in Figure \ref{examp}. 
Embed the graph. The subgraph induced by $\{a,b,1,2,3,4\}$ is $K_6$ and hence has a pair of 3-cycles with non-zero linking number. 

If one 3-cycle is of the form $abx$, where $x \in \{1,2,3,4\}$, then that 3-cycle is consistently oriented. The other component of the link is either $234$ (which is consistent) or of the form $1yz$ with $yz \in \{2,3,4\}$.  If $abx$ has non-zero linking number with $1yz$, note that we may apply Lemma \ref{fix_orient_cor} with $X=abx$, $C_1=1yz$, $C_2=zv1$ and $C_3=1yzv$.   Thus, we have a consistently oriented non-split link.

So, we may assume the vertices $a$ and $b$ are contained in different 3-cycles in the link. Without loss of generality, we may take the pair of linked 3-cycles to be either $a12$ and $b34$ or $b12$ and $a34$. 

In the first case, may apply Lemma \ref{fix_orient_cor} with $X=b34$, $C_1=a12$, $C_2=aw1$ and $C_3=aw12$ to find a consistently oriented cycle $C$ that has nonzero linking number with $b34$.  We may then apply Lemma \ref{fix_orient_cor} with $X=C$, $C_1=b34$, $C_2=b4v$ and $C_3=b34v$ to find a consistently oriented cycle $C'$ that has nonzero linking number with $C$.   By construction, $C$ and $C'$ are disjoint, so we have a consistently oriented nonsplit link.

In the second case, where $b12$ is linked with $a34$, may apply Lemma \ref{fix_orient_cor} with $X=a34$, $C_1=b12$, $C_2=b2v$ and $C_3=b12v$ to find a consistently oriented cycle $C$ that has nonzero linking number with $a34$.  We may then apply Lemma \ref{fix_orient_cor} with $X=C$, $C_1=a34$, $C_2=aw3$ and $C_3=aw34$ to find a consistently oriented cycle $C'$ that has nonzero linking number with $C$.   By construction, $C$ and $C'$ are disjoint, so we have a consistently oriented nonsplit link.

Thus, any embedding of $G$ must contain a consistently oriented nonsplit link. 
We may extend $G$ to a tournament on 8 vertices by adding 3 edges: $bw, vw$ and $av$. Various choices lead to up to eight different directed graphs. Associated score sequences include $(2,3,3,3,4,4,4,5)$, $(2,3,3,3,3,4,5,5)$, $(3,3,3,3,3,4,4,5)$, as well as $(2,3,3,3,3,4,4,6)$. The first three are already known to have intrinsically linked representatives by Proposition \ref{K332_prop}. The last one arises from the orientations $wb$, $wv$ and $va$. 

To obtain the other three score sequences, form $G'$ by modifying the directed graph $G$ from Figure \ref{examp} by changing the orientation of edge $ab$ to be from $b$ to $a$, and adding the  edge $wb$ oriented from $w$ to $b$.

In any embedding of $G'$, we may start with the pair of 3-cycles with non-zero linking number in the subgraph induced by $\{a,b,1,2,3,4\}$ and apply Lemma \ref{fix_orient_cor} as needed to find a consistently oriented link as before.  Thus $G'$ is intrinsically linked as a directed graph.

We may extend $G'$ to a tournament on 8 vertices by adding 2 edges: $wv$ and $va$.   The sequences above correspond to the orientation choices $(wv, va)$, $(vw, va)$ and $(wv, av)$ respectively.
\end{proof}

We now turn our attention to finding linkless score sequences.

\begin{lemma}
Let $T$ be a tournament on 8 vertices with a vertex $v$ of out degree 1 and a vertex $w$ of out degree 6.   If edge $wv$ is oriented from
$w$ to $v$ then $T$ is not intrinsically linked as a directed graph.
\label{lemma_16}
\end{lemma}

\begin{proof}
We may conduct a consistent edge contraction on the edge $vx$ that is oriented from $v$ to $x$ to form a directed graph
$G'$ on 7 vertices, as $v$ is a sink in  $G \setminus vx$.   The vertex $w$ has in degree 1 in $G'$,
so we may conduct a consistent edge contraction on the edge $wx'$ that is oriented from $x'$ to vertex $w$ to form a directed graph
$G''$ on 6 vertices, as $w$ is a source in $G' \setminus x'w$.

Denote the complete symmetric graph $DK_n$.  Any symmetric directed edges in $G''$ must be incident on either the vertex from $vx$ or $x'w$.  Thus $G''$
is a subgraph of $DK_6 \setminus K_4$.  Consider a vertex $y$ in the $K_4$.  It has degree 3, so $y$ must have in degree at least 2 or
out degree at least 2.  Thus, there is a vertex $y$ in $G''$ such that to form $G''$ by deleting edges from $DK_6$, we must delete at least two edges incident on $y$ that are not consistently oriented. Thus, $G''$ is not intrinsically linked by Theorem \ref{FHR_edge_pair}.  As $G''$ is obtained from $T$ by a sequence of consistent edge contractions, $T$ is not intrinsically linked by Theorem \ref{consist_edge_thm}.
\end{proof}

\begin{lemma}  Let $T$ be a tournament on $8$ vertices with score sequence  $S$.  If $S$ contains 0 (or 7), or the subsequence 1,1 (or 6,6), or the subsequence 1, \ldots 5, 5, 6 (or 1, 2, 2 \ldots 6) then $T$ is not intrinsically linked as a directed graph.
\label{lemma_0_reduce}
\end{lemma}

\begin{proof}
Suppose $S$ contains 0 (or 7). Let $v$ be the vertex with out degree 7 (or 0).  A consistently oriented cycle cannot contain $v$.
Thus, if $T$ is intrinsically linked, then
$T \setminus v$ must be intrinsically linked.   However, $T \setminus v$ is a tournament on 7 vertices,
and by Theorem \ref{7_tourneys}, a tournament on 7 vertices is not intrinsically linked as a directed
graph.   Thus, $T$ is not intrinsically linked.

Suppose $S$ contains 1,1 (or 6,6).
By Lemma \ref{1_1_lemma}, $S$ has an intrinsically linked representative if and only if the length 7 score sequence $S'$ does.  But no 7 vertex tournament is intrinsically linked, so $T$ is not intrinsically linked as before.

Suppose $S$ contains 1, 2, 2 \ldots 6.  Label the vertex of out degree 1 as $v$ and the vertex of out degree 6 as $w$.  If the edge between $v$ and $w$ is oriented from $w$ to $v$, then $T$ 
is not intrinsically linked by Lemma \ref{lemma_16}. Thus, we may assume that the edge between $v$ and $w$ is oriented from $v$ to $w$. Conduct a consistent edge contraction on this edge
and label the resulting vertex $7$.

The resulting directed graph $G$ has 7 vertices, and has symmetric edges only incident on vertex $7$.  Thus, $G \setminus 7$ is a tournament $T'$ on 
6 vertices.
All other vertices $x$ of $T$ have edges directed from $x$ to $v$ and from $w$ to $x$. Thus, the score sequence of $T'$ contains $( \dots 4,4)$.
By Lemma \ref{sub_tourney_0_11}, $G$ is not intrinsically linked as a directed graph. As $G$ was obtained from $T$ by a consistent edge contraction, Theorem \ref{consist_edge_thm} implies $T$ is not intrinsically linked as a directed graph.

\end{proof}

\begin{prop}
Let $T$ be a tournament on 8 vertices such that the sum of the out degree of four of the vertices is 8.  Then $T$ is not intrinsically linked as a directed graph.
\label{2_out_edges_prop}
\end{prop}

\begin{proof}

Label the vertices whose out degree sums to 8 as $\{A,B,C,D\}$.  Label the other vertices $\{a,b,c,d\}$.  The vertices $\{A,B,C,D\}$ have a total out degree of 8.  As these vertices form a $K_4$, there are exactly 2 edges $e_1$ and $e_2$ oriented from a vertex in $\{A,B,C,D\}$ to a vertex in $\{a,b,c,d\}$.

Suppose that $e_1$ and $e_2$ have a common vertex.  Each set of vertices $\{a,b,c,d\}$ and $\{A,B,C,D\}$ form a $K_4$.   We may choose an embedding $f$ of $T$ such that these $K_4$ graphs are contained in disjoint 2-spheres and are paneled.   That is, each 3-cycle in $f(\{a,b,c,d\})$ and $f(\{A,B,C,D\})$ bounds a disk.  Then if $L= L_1 \cup L_2$ is a link in $f(T)$ where $L_1$ and $L_2$ are consistently oriented, and $L_1$ or $L_2$ contains vertices from only one of the partitions, then $L$ is a split link.  Thus if $L$ is a non-split consistently oriented link, each $L_i$ must contain at least one vertex from both $\{a,b,c,d\}$ and $\{A,B,C,D\}$.  Suppose that $e_1 \in L_1$.  Then $e_2 \notin L_2$, as $L_1$ and $L_2$ are disjoint.  So $L_2$ is not consistently oriented, a contradiction.  

Thus, we may assume that $e_1$ and $e_2$ do not have a common vertex.  Embed $T$ as shown in Figure \ref{2_out_edges_fig}. We will show that this embedding contains no non-split consistently oriented links.

\begin{figure}
\includegraphics[scale=.6]{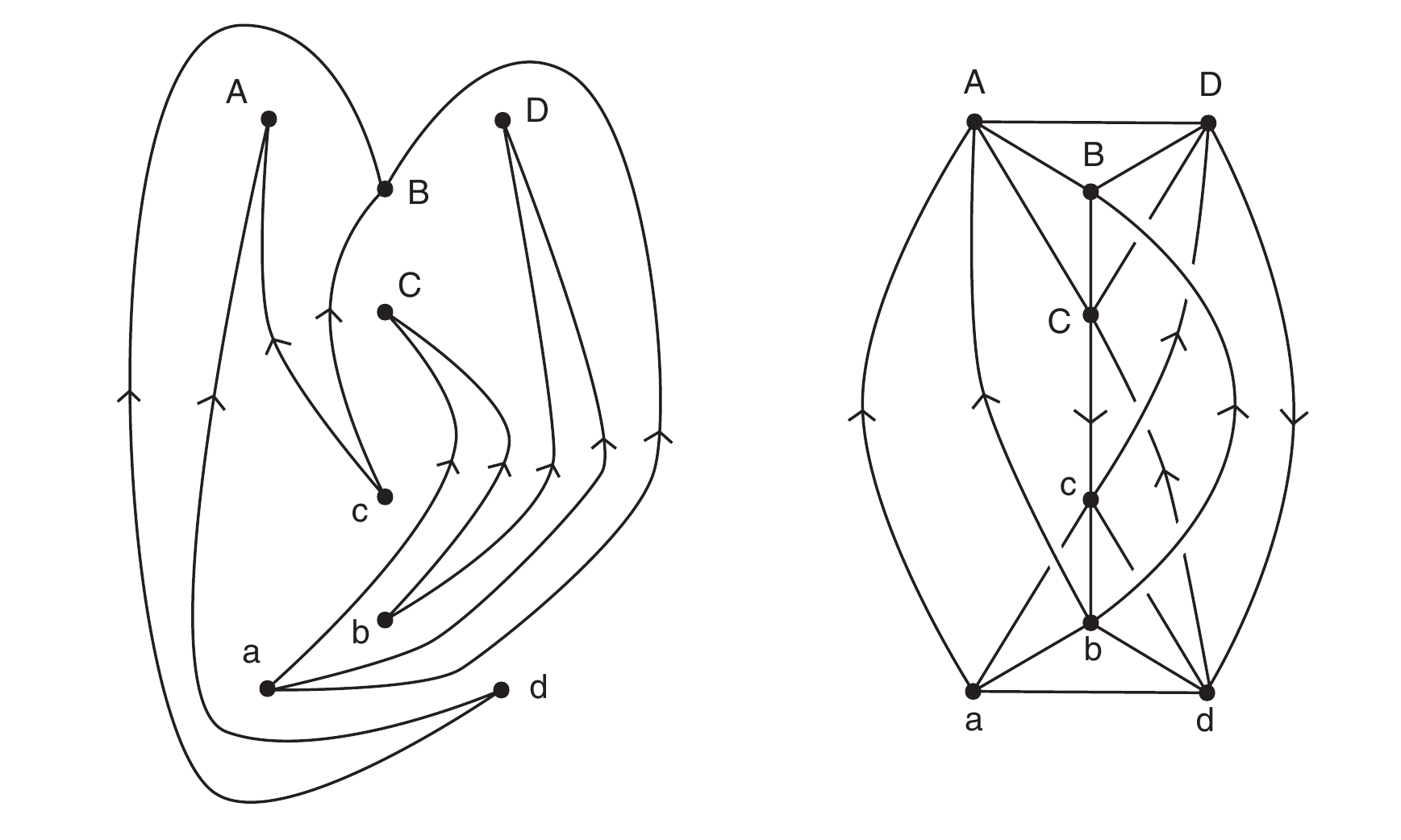}
\caption{An embedding of $T$.  The vertices lie in the plane $z=0$.  The edges shown on the left side of the figure have $z<0$ and the edges shown at right have $z>0$.}
\label{2_out_edges_fig}
\end{figure}

Suppose that $L = L_1 \cup L_2$ is a consistently oriented link where $L_1$ and $L_2$ are 4-cycles.  Suppose $L_1 \subset \{A,B,C,D\}$ or $L_1 \subset \{a,b,c,d\}$ then $L_1$ bounds a disk, and $L$ is split. Similarly for $L_2$. Suppose $L_1$ has one vertex from one partition and three from the other.   As $L_1$ is consistently oriented, it must contain edge $Cc$ or $Dd$.  Further, if $L_1$ contains $Cc$, then $L_2$ contains $Dd$ as $L_2$ is consistently oriented.  Hence $L_1$ must be of the form $ABCc, BACc, ABDd, BADd, abcC, bacC, badD$ or $abdD$.  Note that the 3-cycles $ABC, ABD, abc,$ and $abd$ bound disks, and hence the 4-cycles $ABCc, BACc, ABDd, BADd, abcC, bacC, badD$ and $abdD$ are isotopic to the 3-cycles $ACc, BCc, ADd, BDd, acC, bcC, bdD$ and $adD$ respectively.   We will show that each of these 3-cycles bounds a disk below, when we consider potential links that include a 3-cycle.

So we may assume that $L_1$ has two vertices from each partition, and that $L_2$ does as well.  As $L_1$ is consistently oriented, it must contain edge $Cc$ or $Dd$, and $L_2$ must contain the other.  Thus, the link $L$ is one of: $DdbA,  CcaB$ or $DdbB, CcaA$ or $DdaA, CcbB$ or $DdaB, CcbA$.  Examining Figure \ref{2_out_edges_fig}, we see that each of these links is split.

Thus, we may assume that $L_1$ is a 3-cycle.

If $L_1 \subset \{A,B,C,D\}$ or $L_1 \subset \{a,b,c,d\}$, then either $L_1$ bounds a disk, or $L_1$ is one of $acd$ or $ACD$.  If $L_2$ is of the form $bXY, bXYZ, bXYZW, Bxy, Bxyz$ or $Bxyzw$ then $L_2$ is not consistently oriented.  Thus $L_2\subset \{A,B,C,D\}$ or $L_2 \subset \{a,b,c,d\}$, and $L$ is split.

So we may assume that $L_1$ is of the form $cCX, cCx, dDX,$ or $dDx$.  The twelve possible cycles $cCA, cCB, cCD, acC, bcC, dcC, dDA, dDB, adD, bdD,$ and $cdD$ all bound disks.
\end{proof}

We have the immediate corollary.

\begin{cor}
Let $T$ be a tournament with score sequence $(1,2,2,3,5,5,5,5)$, $(2,2,2,2,5,5,5,5)$ or $(2,2,2,2,4,5,5,6)$.   Then $T$ is not intrinsically linked as a directed graph.
\label{cor_8_out}
\end{cor}

\begin{prop}
Let $T$ be a tournament on $8$ vertices with score sequence \linebreak $(1,2,3,3,4,5,5,5)$, $(1,3,3,3,4,4,4,6)$, $(1,3,4,4,4,4,4,4)$, $(1,2,3,3,4,4,5,6)$, \linebreak
$(1,2,2,4,4,5,5,5)$, 
 $(1,2,4,4,4,4,4,5)$, $(1,3,3,3,3,5,5,5)$, or a score sequence dual to one of these. Then $T$ is not intrinsically linked as a directed graph.
\label{prop_special_linkless}
\end{prop}

\begin{proof}

Suppose $T$ has score sequence $(1,2,3,3,4,5,5,5)$.
Label the vertex of out degree 1 as $x$, the vertex of out degree 2 as $y$, the vertices of out of degree 3 as $z_1, z_2$, the vertex of out degree 4 as $w$, and the vertices out degree 5 as $v_1, v_2, v_3$.

There is a single edge oriented from $x$ to a vertex $u$, and hence we may conduct consistent edge contraction of edge $xu$ to form a 7 vertex directed graph $G$.  Label the newly formed vertex $7$. Deleting 7 results in a tournament $T'$ on 6 vertices.   We can examine the score sequence $S'$ of $T'$ under the possible choices of $u$ to show that $G$ has a linkless embedding.

Suppose $u = y$.  
Vertex $y$ has out degree 2, so there are two edges $ys$ and $yt$ oriented from $y$ to $s$ and $t$. Suppose $z_i \neq s$, then we have
edge orientations $z_1x, z_1y, z_2x, z_2y$, so $1,1 \in S'$, and $G$ is linkless by Lemma \ref{sub_tourney_0_11}.  So we may assume $s=z_1$.
Suppose $t=z_2$, then we have edge orientations $z_1x, yz_1, z_2x, yz_2$, so $S' = (2,2,2,3,3,3)$ and $G$ is linkless by Lemma \ref{sub_tourney_222333}.
Suppose $t=v_1$, then we have edge orientations $z_1x, yz_1, v_1x, yv_1$, so $S' = (1,2,2,3,3,4)$ and $G$ is linkless by Lemma \ref{sub_tourney_222333}.
The only remaining possibility (up to symmetry) is $s=z_1, t=w$. In this case we have edge orientations 
$z_1x, yz_1 wx, yw$, and $S' = (1,2,3,3,3,3)$.  By Lemma \ref{sub_tourney_123333} $G$ is linkless  unless $T'$ has edge orientations $ba, ac, da,$ $ea, fa, bc,$ $db, eb, fb,$ $cd, ce, cf,$ $de, fd, ef$.

We have edge orientations $v_1x$ and $v_1y$.  Thus there are two edges directed from $v_1$ to $7$ in $G$, and as they have the same orientation, we may delete one of them.  As $v_1$ has outdegree 3 in $T'$, we may label it $d$. Thus, there is a single edge directed from $d$ to $7$ in $G$, and so $G$ is linkless by Lemma \ref{sub_tourney_123333}.

So $u$ must be one of $w, z_1,$ or $v_1$, and we may use similar arguments to show $G$ has a linkless embedding in these cases as well. Since we may find a linkless embedding for $G$ no matter the structure of $T$, and $G$ is obtained from $T$ by consistent edge contraction,
$T$ is not intrinsically linked by Theorem \ref{consist_edge_thm}.

We may make similar arguments for the other score sequences.  
\end{proof}

There are 167 score sequences for 8 vertex tournaments, and we are able to classify all but 5 of them, as shown in Table \ref{8_vert_table}.  The remaining unclassified score sequences are shown in Remark \ref{unkown_8_verts} below. We conjecture these remaining score sequences to be linkless.

\begin{remark}
\label{unkown_8_verts}
The results of this section allow the classification of all score sequences for 8 vertex tournaments except those below, shown in dual pairs.
\\
$
(2, 2, 2, 3, 4, 5, 5, 5) (symmetric)\\
(2, 2, 2, 4, 4, 4, 5, 5) (2, 2, 3, 3, 3, 5, 5, 5)\\
(2, 2, 4, 4, 4, 4, 4, 4) (3, 3, 3, 3, 3, 3, 5, 5)\\
$
\end{remark}

\section{Tournaments on 9 or more vertices}
\label{9_vert_section}

The results of Section \ref{8_vert_section} allow us to classify most score sequences for 8 vertex tournaments.   In this section we will attack the classification of score sequences for larger tournaments in two ways.   First, we will develop some techniques for reducing the classification of a length $n$ score sequence to classifying a length $n-1$ score sequence, which allows us to apply the classification of score sequences for 8 vertex tournaments to the larger cases.   Second, we will extend some of the techniques from Section \ref{8_vert_section} to longer score sequences directly.

With these techniques we can classify all but 37 of the 490 score sequences for 9 vertex tournaments, all but 150 of the 1486 score sequences for 10 vertex tournaments, and all but 512 of the 4639 score sequences for 11 vertex tournaments.

We first look at methods for classifying score sequences for $n$ vertex tournaments based on score sequences for smaller tournaments.

\begin{lemma}
If an $n$ vertex tournament $T$ has a score sequence that contains 0 or $n-1$, that tournament has an intrinsically linked
representative if and only if $T' = T\setminus v$ does, where $T'$ is a tournament on $n-1$ vertices and $v$ is a vertex of out degree 0 or $n-1$.
\label{0_n_lemma}
\end{lemma}

\begin{proof}
As $v$ has out degree 0 or $n-1$, it cannot be contained in a consistently oriented cycle. Thus, $T$ is intrinsically linked as a directed graph if and only if $T \setminus v$ is.
\end{proof}

\begin{obs}
If an $n-1$ vertex tournament $T'$ is intrinsically linked, then an $n$ vertex tournament $T$ formed by adding a vertex of any out degree to $T'$ is intrinsically linked.
\label{linked_extend}
\end{obs}

Observation \ref{linked_extend} implies that given a score sequence for an 8 vertex tournament that has an intrinsically linked representative, we can form many score sequences for 9 vertex tournaments that have intrinsically linked representatives.   For example, (3,3,3,3,4,4,4,4) has an intrinsically linked representative.   We may add a vertex of out degree 7 to show that the 9 vertex score sequences (3,3,3,3,4,4,4,5,7) and (3,3,3,4,4,4,4,4,7) have intrinsically linked representatives.  This allows us to classify many score sequences for larger tournaments.

We may also make use of consistent edge contraction.

\begin{lemma}
The $n$ vertex score sequence $(1,1, s_3,s_4, \ldots s_n)$ has an intrinsically linked representative if and only if the $n-1$ vertex score sequence $(1, s_3-1, s_4-1, \ldots s_n-1)$ does.  Similarly, the $n$ vertex score sequence $(s_1, s_2, \ldots n-2, n-2)$ has an intrinsically linked representative if and only if the $n-1$ vertex score sequence $(s_1, s_2, \ldots n-3)$ does.
\label{1_1_lemma}
\end{lemma}

\begin{proof}
Let $T$ be an $n$ vertex tournament with score sequence $(1,1, s_3,s_4, \ldots s_n)$. Let $v$ and $w$ be the vertices of out degree 1.
We may assume that the edge  $vw$ is oriented from $v$ to $w$, and we may conduct a consistent edge contraction on this edge, as $v$ is
a sink in $T \setminus vw$. Let $wx$ be the unique edge directed from $w$ to another vertex. Let $G$ be the directed graph obtained from
$T$ by consistent edge contraction on edge $vw$. We will abuse notation and call the resulting vertex $w$.

$T$ is intrinsically linked if and only if $G$ is, by Theorem \ref{consist_edge_thm}.  The digraph $G$ has a single pair of symmetric directed edges between
$w$ and $x$.  Suppose that the edge $xw$ is used in a consistently oriented cycle $c$. As $w$ has outdegree 1, $wx \in c$.   Thus $c$ is the 2-cycle $wx-xw$.
As we may choose to embed $G$ so that $wx-xw$ bounds a disk, $G$ is intrinsically linked if and only if $T' = G \setminus xw$ is. As $T'$ is tournament on $n-1$ vertices, with score sequence $(1, s_3-1, s_4-1, \ldots s_n-1)$, we have the result.

We may apply the same argument with all of the edge orientations reversed to show the result for $(s_1, s_2, \ldots n-2, n-2)$ and $(s_1, s_2, \ldots n-3)$
\end{proof}

We are able to generalize Proposition \ref{2_out_edges_prop} to find examples of linkless score sequences for tournaments with up to 14 vertices.

\begin{lemma}
Let $T$ be an $n$ vertex tournament.  If the first $m$ values of the score sequence sum to ${m\choose2}$ and $m < 8$ and $(n-m) < 8$, then $T$ is not intrinsically linked as a directed graph.
\label{linkless_split}
\end{lemma}

\begin{proof}
Denote the first $m$ vertices $M = \{ v_1, v_2 \dots v_m\}$. The vertices of $M$ have total out degree ${m\choose2}$, and hence all edges between $T' = T \setminus M$ and $M$ must be oriented from $T'$ to $M$.   Thus no consistently oriented cycle may contain a vertex from $M$ and a vertex from $T'$.

As $m<8$, by Theorem \ref{7_tourneys} we may embed $M$ so that there is no consistently oriented non-split link using only vertices of $M$.  Enclose this embedding in the interior of an $S^2$.

Choose a second $S^2$ disjoint from the one containing $M$, and embed $T'$ inside it.  As $(n-m)<8$, we may choose an embedding of $T'$ so that there is no consistently oriented non-split link using only vertices of $T'$.  

We may now extend this to an embedding of $T$ by embedding the edges between $T'$ and $M$ arbitrarily.

If $L_1, L_2$ is a pair of consistently oriented cycles in this embedding of
$T$, then either $L_1, L_2$ are both in $T'$ or both in $M$, and hence split, or $L_1$ is in $M$ and $L_2$ is in $T'$.   These embeddings are divided by an $S^2$, so this link is split as well.
\end{proof}

We are able to strengthen Lemma \ref{linkless_split} in the special case of 9 and 10 vertex tournaments.

\begin{lemma}
Let $T$ be a tournament.  If $T$ has 9 vertices and the first four elements of its score sequence sum to 7, or $T$ has 10 vertices and the first five elements of its score sequence sum to 11, then $T$ is not intrinsically linked as a directed graph.
\label{nine_vert_7_out}
\end{lemma}

\begin{proof}
Suppose $T$ is a tournament on 10 vertices, and that we may partition the vertices as $A = \{1,2,3,4,5\}$ and $B=\{1',2',3',4',5'\}$, where the total out degree of the vertices in $A$ sums to 11.  

As the there are 10 edges between vertices of $A$, and these vertices have total out degree 11, there is exactly 1 edge directed from $A$ to $B$.  Chose the vertex labels so that this is edge $22'$.

Form an embedding of $T$ as follows.
The vertices of $B$ and the edges between them form $K_5$.  Embed this $K_5$ so that there is a single crossing, edge $2'4'$ crossing over edge $1'3'$, and all of $B$ lies in the plane $z=0$ except edge $2'4'$ which lies slightly above it. Embed $A$ so that it is a reflection of $B$ across $z=\frac{1}{2}$.  
Embed the remaining edges of $T$ arbitrarily between the planes $z=0$ and $z=1$.  

Let $L = L_1 \cup L_2$ be a consistently oriented link.  Note that if $L_1 \subset A$ then $L_2 \not \subset A$.  Similarly for $B$.
Suppose $L_1 \subset A$ and $L_2 \subset B$.  Then $L$ is split.

There is exactly one edge from $A$ to $B$.  If both $L_1$ and $L_2$ contain vertices from both $A$ and $B$, then at least one of them must have an inconsistent orientation.

Thus we may assume that $22' \in L_1$, and as $L_1$ and $L_2$ are disjoint, $L_2$ cannot contain $24$ or $2'4'$.  If $L_2 \subset B$ (or $L_2 \subset A$) then after perturbing $L_1$, it lies above (or below) $L_2$, and hence $L$ is split.

The proof of the 9 vertex case is similar, with $A$ having 4 vertices and a planar embedding.
\end{proof}

\section{Further Results}
\label{further_sec}

In Section \ref{8_vert_section}, we found a large number of linkless score sequences for 8 vertex tournaments.  Lemma \ref{0_n_lemma} implies that these may be extended to form linkless score sequences for tournaments of any order by iteratively adding vertices of out degree 0.   In particular, we have the following bound on the minimum number of linkless score sequences for any $n \geq 8$.

\begin{prop}
Let $k$ be the number of linkless score sequences for 8 vertex tournaments, and $n \geq 8$. Then for $n$ vertex tournaments, there are at least $(n-7)(k-1)+1$ linkless score sequences.
\label{always_linkless}
\end{prop}

\begin{proof}
Given a linkless score sequence $S=(s_1, s_2, s_3, s_4, s_5, s_6, s_7, s_8)$, we may construct a linkless score sequence on $n$ vertices by iteratively adding $m$ vertices of degree 0 and $n-8-m$ vertices of maximal degree.

Adding a vertex of degree 0 to $S' = (s_1, s_2, \ldots  s_{m-1}, s_m)$ gives $S'' = (0, s_1+1, s_2+1, \ldots  s_{m-1}+1, s_{m}+1)$, and adding a vertex of degree $m$ gives $S''' = (s_1, s_2, \ldots s_{m-1}, s_m, m)$. Hence,  $S'' = S'''$ if and only if $s_1 = 0$, $s_2 = s_1+1 \ldots s_m+1 = m$.   This occurs only when $S$ is the score sequence of a transitive tournament.  

Note that the order in which the new vertices are added does not matter to the final score sequence, only the number $m$ added of minimal degree.  Thus there are $n-7$ possible results for each sequence $S$, corresponding to adding between $0$ and $n-8$ vertices of minimal degree.  

Excluding the transitive tournament, we have $k-1$ linkless score sequences for 8 vertex tournaments, each of which can be extended to $n-7$ linkless score sequences for $n$ vertex tournaments.  Each of these $(n-7)(k-1)$ sequences is unique, as they contain a subsequence of the form $( \ldots s_1+m, s_2+m, s_3+m, s_4+m, s_5+m, s_6+m, s_7+m, s_8+m, \ldots)$, where $(s_1, s_2, s_3, s_4, s_5, s_6, s_7, s_8)$ a linkless score sequence for 8 vertex tournaments, and $m$ is the number of vertices of minimal degree that were added in constructing $S$.

Including the transitive tournament on $n$ vertices gives the bound.
\end{proof}

By the classification in Section \ref{8_vert_section}, $147 \leq k \leq 152$ for Proposition \ref{always_linkless}.   A careful application of Lemma \ref{1_1_lemma} may improve the lower bound in Proposition \ref{always_linkless}.

While Proposition \ref{always_linkless} implies that for tournaments of order $n$ there are at least $O(n)$ linkless score sequences, we conjecture that as the tournaments become large, score sequences that guarantee that a tournament is not intrinsically linked as a directed graph become vanishingly rare.  More precisely, let $u(n)$ be the number of linkless score sequences for tournaments on $n$ vertices and let $s(n)$ be the total number of score sequences  for tournaments on $n$ vertices.   We conjecture that $\frac{u(n)}{s(n)}$ goes to 0 as $n$ goes to infinity.  

As support for this conjecture, note that the number of score sequences of length $n$ is known to grow as  $O(4^nn^{-\frac{5}{2}})$ \cite{wk} \cite{kp}, and Observation \ref{linked_extend} implies that the number of score sequences with intrinsically linked representatives should grow rapidly as well.  In fact, for 8 vertex score sequences, we have at least 15 score sequences with intrinsically linked representatives, at least 131 for 9 vertex score sequences, at least 660 for 10 vertex score sequences, and at least 2719 for 11 vertex score sequences.  

\medskip

Given a tournament $T$ with a linkless score sequence, we know that $T$ must have a linkless embedding.  When a tournament has a score sequence with an intrinsically linked representative, the intrinsic linking status of the tournament is unclear.  Specifically, let $T$ be a tournament with score sequence $S$, where $S$ has an intrinsically linked representative.  Then $T$ may be intrinsically linked as a directed graph, or $T$ may have an embedding with no consistently oriented nonsplit link, as the following proposition shows.

\begin{prop}
Let $T$ be a tournament on 8 vertices with score sequence \linebreak $(1,3,3,3,3,4,5,6)$.   Then $T$ may be intrinsically linked as a directed graph, or $T$ may have an embedding with no consistently oriented nonsplit link.
\end{prop}

\begin{proof}
By Proposition \ref{flex_prop}, there exists a tournament $T$ with score sequence \linebreak $(1,3,3,3,3,4,5,6)$ where $T$ is intrinsically linked as a directed graph.

Let $T'$ be a tournament with score sequence $(1,3,3,3,3,4,5,6)$.  Let $v$ be the vertex of out degree 1, and $w$ the vertex of out degree 6, and let edge $wv$ be oriented from $w$ to $v$.  Then $T'$ is not intrinsically linked as a directed graph by Lemma \ref{lemma_16}.
\end{proof}

Thus, in this work we have divided score sequences into two classes--those that imply the tournament must have a linkless embedding, and those that leave the intrinsic linking status of the tournament uncertain.  We conjecture the existence of a third class--specifically, we conjecture that there exists a score sequence $S$ such that any tournament with score sequence $S$ must be intrinsically linked as a directed graph.  What is the smallest number of vertices required for such a score sequence?


\begin{thebibliography}{99}

\bibitem{REU} S. Chan, A. Dochtermann, J. Foisy, J. Hespen, E. Kunz, T. Lalonde, Q. Loney, K. Sharrow, N. Thomas, \emph{Graphs with disjoint links in every spatial embedding}, J. Knot Theory Ramifications \textbf{13} No. 6 (2004)  737-748.

\bibitem{cg} J. H. Conway, C. McA. Gordon, \emph{Knots and links in spatial graphs}, J. Graph Th. \textbf{7} (1983) 446-453

\bibitem{flemfoisy} T. Fleming, J. Foisy, \emph{Intrinsically knotted and 4-linked directed graphs}, J. Knot Theory Ramifications \textbf{27} No. 6 (2018) 1850037-1 to 1850037-18

\bibitem{flemfoisy2} T. Fleming, J. Foisy, \emph{Intrinsic linking and knotting in tournaments},  J. Knot Theory Ramifications 28  (2019) no. 12, 1950076-1 to 1950076-18 

\bibitem{fmellor} T. Fleming, B. Mellor, \emph{Counting links in complete graphs}, Osaka J. Math. \textbf{46} (2009), 173-201

\bibitem{FHR} J. Foisy, H. Howards, N. Rich \emph{Intrinsic linking in directed graphs}, Osaka J. Math. \textbf{52} (2015) 817-831

\bibitem{kp} J. H. Kim, B. Pittel, \emph{Confirming the Kleitman-Winston conjecture on the largest coefficient in a q-Catalan number}, J. Comb. Theory A \textbf{92} (2000) 197-206

\bibitem{landau} H. G. Landau, \emph{On dominance relations and the structure of animal societies. III. The condition for a score structure}, Bull. Math. Biophysics, \textbf{15} (1953) 143-148

\bibitem{mader} W. Mader, \emph{Homomorphieeigenschaften und mittlere Kantendichte von Graphen}, Mathematische Annalen \textbf{174} No. 4 (1967) 265-268

\bibitem{mnp} T. Mattman, R. Naimi, B. Pagano, \emph{Intrinsic linking and knotting are arbitrarily complex in directed graphs}, arXiv:1901.01212

\bibitem{nt} J. Ne\v{s}et\v{r}il, R. Thomas, \emph{A note on spatial representation of graphs} Commentat. Math. Univ. Carolinae \textbf{26} No. 4 (1985) 655-659

\bibitem{rs} N. Robertson,  P. Seymour, \emph{Graph minors. XX. Wagner's conjecture} J. Comb. Theory B \textbf{92} (2004) 325-357

\bibitem{rst} N. Robertson, P. Seymour, R. Thomas, \emph{Sachs' linkless embedding conjecture}, J. Comb. Theory  B \textbf{64} (1995) 185-277

\bibitem{sachs} H. Sachs, \emph{On spatial representations of finite graphs}, in: A. Hajnal, L. Lovasz, V.T.  S\'os (Eds.), Colloq. Math. Soc. J\'anos Bolyai, Vol. 37, North-Holland, Amsterdam, (1984) 649-662

\bibitem{wk} K. J. Winston and D. J. Kleitman, \emph{On the asymptotic number of tournament score sequences}, J. Comb. Theory A \textbf{35} (1983) 208-230

\end{thebibliography}
\end{document}